\documentclass[11pt]{article}
\usepackage{amsmath,amsthm}
\usepackage{graphicx,psfrag,epsf}
\usepackage{enumerate}
\usepackage[longnamesfirst]{natbib}

\usepackage{color,amssymb}
\usepackage{epsfig}
\usepackage{rotating}
\usepackage{booktabs}
\usepackage{bm}
\usepackage{setspace}
\usepackage{dcolumn}
\usepackage{authblk}

\newcolumntype{d}{D{.}{.}{-1}}

\theoremstyle{plain}

\newcommand{\be}{\begin{eqnarray}}
\newcommand{\ee}{\end{eqnarray}}
\newcommand{\bq}{\begin{eqnarray*}}
	\newcommand{\eq}{\end{eqnarray*}}

\newcommand{\I}{\mathbb{I}}

\DeclareMathOperator{\sgn}{sgn}
\DeclareMathOperator{\Prob}{P}

\title{\bf Fourier-type monitoring procedures \\ for strict stationarity}
\author[1]{Sangyeol Lee}
\author[2]{Simos G. Meintanis\thanks{On sabbatical leave from the National and Kapodistrian University of Athens, Athens, Greece.}}
\author[2,3]{Charl Pretorius}
\affil[1]{Department of Statistics, Seoul National University, Seoul, South~Korea}
\affil[2]{Unit for Business Mathematics and Informatics, North-West University, Potchefstroom, South~Africa}
\affil[3]{Department of Probability and Mathematical Statistics, Charles University, Prague, Czech~Republic}
\date{\today}

\begin{document}
	
	\maketitle
	
	\bigskip
	\begin{abstract}
		We consider model-free monitoring procedures for strict stationarity of a given time series.
		The new criteria are formulated as L2-type statistics incorporating the empirical characteristic
		function. Asymptotic as well as Monte Carlo results are presented. The new methods are also employed
		in order to test for possible stationarity breaks in time-series data from the financial sector.
	\end{abstract}
	
	
	\section{Introduction}\label{sec_1}
	
	The notion of stationarity plays a very
	important role in statistical modeling. In its weakest form of first or second order stationarity it
	implies that the mean or second moment, respectively, are time invariant; see for instance
	\citet{xiao07}, \citet{dwivedi11}, and \citet{SSRAO15}.  A more general
	related notion is that of $p^{{\rm{th}}}$ order (weak) stationarity which
	requires that all joint product moments of order (up to) $p$ are time
	invariant.  Most studies of stationarity are restricted to some form
	of weak stationarity, which of course is the most suitable concept
	for linear time series. On the other hand, the
	property of strict stationarity states that not only moments, but
	the entire probabilistic structure of a given series is time
	invariant. This property is of great relevance with non-linear
	time series in which low order moments are not sufficient for the dynamics of
	the series, not to mention the case of heavy-tailed time series lacking higher moments \citep{loretan94,andrews09}.
	Another divide is between parametric and non-parametric tests for stationarity, with the first class
	containing the majority of procedures in earlier literature. There is also the methodological approach
	that categorizes methods which operate either in the time or in the frequency domain. As the existing
	literature is vast, we provide below only a selected set of references which is by no means exhaustive.
	
	In econometrics the majority of earlier tests for weak stationarity are either tests for stationarity against
	unit root, or alternatively, tests of a unit root against  stationarity, with the KPSS and the Dickey-Fuller
	tests being the by far most popular ones and having enjoyed a number of generalizations;
	see for instance \citet{gilalana03}, \citet{giraitis03}, \citet{mueller05}, and \citet{horvath14}. In fact these tests, although developed
	with a particular alternative in mind, they have often been employed more generally to explore structural change.
	When it comes to testing for strict stationarity in parametric time series there exist many tests. These tests typically reduce
	to testing for a given portion of the parameter space and may often readily be extended to monitoring procedures. We indicatively mention
	here the works for ARMA, GARCH, and DAR models by \citet{bai94}, \citet{francq12} and \citet{guo16}, respectively. On the other
	hand testing methods for strict stationarity are scarce when performed within a purely nonparametric framework. It appears that \citet{kapetanios09} was
	the first to address the issue of testing strict stationarity of the marginal distribution of an arbitrary time series. His work employs smoothing
	techniques in order to estimate the density and a bootstrap implementation of the procedure. There is also the method of \citet{hong2016} which is based on the joint characteristic function and the papers of \citet{busetti10} and \citet{lima13}
	which test for constancy (of a discretized version) of the marginal quantile process.
	The interest in testing for stationarity rests on the fact that modelling, predictions and other inferential procedures are invalid if this assumption is violated. However, although strict stationarity is widely assumed in the literature, it
	is not truly a realistic assumption when one observes a given time
	series over a long period of time. On the contrary it is expected
	that institutional changes cause structural breaks in the stochastic
	properties of certain variables, particularly in the macroeconomic
	and financial world. In this connection, monitoring the stationarity of a stochastic process seems to be of
	an even greater importance than testing. In this paper we propose a sequential
	procedure for strict stationarity. Our approach uses the characteristic function (CF) as the main tool. The advantage of using this function is that the CF can be estimated purely non-parametrically without the use of smoothing techniques. Moreover, and unlike the case of estimating the joint distribution function, the estimate of the joint CF is easy to obtain and when viewed as a process it is continuous in its argument. These features offer specific theoretical and computational simplifications which are particularly important in the multivariate context, and are clearly reflected in the competitiveness of the resulting methods over classical ones; see for instance \citet{pinske98}, \citet{su07}, \citet{hlavka12}, and \citet{matteson14}. The asymptotics of the detector statistic involves some weak convergence
	theorems of stochastic integrals to achieve its limiting
	distribution. Here, we establish the weak convergence theorem in a
	more general framework for potential usage in various on-line
	monitoring procedures and retrospective change point tests.
	For more details on the monitoring procedure for the autocovariance
	function of a linear process, we refer to \citet{na11}.
	
	The remainder of the paper is as follows. In Section 2 we introduce
	the basic idea behind the proposed procedures. Section 3 presents
	the corresponding detector statistics and in Section 4 we study the
	large sample behavior of the new methods. As the limit null
	distribution of our detector is complicated, we propose in Section 5
	a resampling procedure in order to actually carry out the suggested
	monitoring method. The results of a Monte Carlo study for the
	finite-sample properties of the methods are presented in Section 6.
	Some real-world empirical applications to financial market data are
	presented and discussed in Section 7. Finally, we end in Section 8
	with conclusions and discussion. The proofs of the lemmas and
	theorems in Section 4 are all provided in the Appendix.
	
	\section{ECF statistics}
	Let $\{X_t\}_{t\in \mathbb N}$ be an arbitrary time series, and
	write $F_{t}(\cdot)$ for the corresponding distribution function
	(DF) of the $m$-dimensional variable
	$\Upsilon_t=(X_{t-(m-1)},...,X_{t})'$, $m\geq1$.  We are interested in the behavior of the distribution of $\Upsilon_t$ over time, i.e. to the monitoring the joint distribution of the observations $X_t$ of a given dimension $m$. The null hypothesis is stated as,
	\begin{equation}\label{NULL} {\cal{H}}_0: F_{t}\equiv
	F_{m} \ {\rm{for \ all }} \ t\geq m+1,\end{equation} against the
	alternative
	\begin{equation}\label{ALT} {\cal{H}}_1: F_{t}\equiv F_{m}, \ t\leq t_0 \text{ and }  F_t(y) \neq F_{m}(y), \ t>t_0,  \end{equation}
	for some $y \in \mathbb R^m$, where $F_m$, $F_t$, as well as the threshold $t_0$ are
	considered unknown. Clearly $m=1$ corresponds to monitoring the marginal distribution of $X_t$, $m=2$ corresponds to the joint bivariate distribution of $(X_{t-1},X_t)^\prime$, and so on. As it is typical in monitoring studies we assume that there exist a set of observations $X_1,...,X_T$ (often termed
	{\it{training data}}) which involve no change, and 
	monitoring starts after time $T$.
	
	To motivate our procedure let $\psi_{Y}(u):=E(e^{iu^\prime Y}), \ i=\sqrt{-1}$, be the characteristic function (CF) of an arbitrary random vector
	$Y$. We will compare a nonparametric estimator of the joint
	CF $\psi_{\Upsilon_j}(u)$, of $\Upsilon_j, \
	j=m, m+1,...,T$, with the the same estimator obtained from
	observations beyond time $T$.
	Then, the quantity of interest is
	\begin{equation} \label{DIST}
	D_{T,t}=\int_{\mathbb R^m}
	\left|\widehat \psi_{T}(u)-\widehat \psi_{T+t}(u)\right|^2 w(u) du,
	\end{equation}
	where
	\begin{equation} \label{ECF}  \widehat \psi_J(u)=\frac{1}{J-m+1} \sum_{j=m}^J e^{i u' \Upsilon_j}, 
	\end{equation}
	is the empirical characteristic function (ECF) computed from the observations
	$\Upsilon_{j}, \ j=m,...,J$, and $w(\cdot)$ is a weight function
	which will be discussed below.
	
	Our motivation for considering (\ref{DIST}) as our main tool is that
	the null hypothesis (\ref{NULL}) may equivalently be stated as
	\begin{equation}\label{NULL1} {\cal{H}}_0: \psi_{\Upsilon_{t}}\equiv \psi_{\Upsilon_m} \ {\rm{for \ all }} \ t\ge m+1,\end{equation}
	and therefore we expect $D_{T,t}$ to be `small' under the null
	hypothesis (\ref{NULL1}).  Moreover, and unlike equivalent
	approaches based on the empirical DF, the ECF approach enjoys the
	important feature of computational simplicity. Specifically, by
	straightforward algebra we have from (\ref{DIST}), \
	\begin{multline} \label{DIST1}
	D_{T,t}=\frac{1}{(T-m+1)^2} \sum_{j,k=m}^T W_{j,k}+\frac{1}{(T-m+1+t)^2}
	\sum_{j,k=m}^{T+t} W_{j,k} \\
	-\frac{2}{(T-m+1)(T-m+1+t)} \sum_{j=m}^T
	\sum_{k=m}^{T+t} W_{j,k},
	\end{multline}
	where $W_{j,k}=W(\Upsilon_j-\Upsilon_k)$ with
	\begin{equation} \label{INT} W(x)=\int_{\mathbb R^m} \cos(u'x) w(u) du, \end{equation}
	and a large (positive) value of the test criterion indicates violation of the null hypothesis.   	
	A standard choice is to set $w(u)=e^{-a\|u\|^2}, \ a>0$, which leads to
	\begin{eqnarray}\label{EXPW}
	W(x)=\left(\frac{\pi}{a}\right)^{m/2} e^{-\|x\|^2/4a},
	\end{eqnarray}
	and hence renders our statistic in closed-form. 
	
	Another interesting
	choice results by considering the statistic $\widetilde
	D_{T,t}=-D_{T,t}$ (in which case of course large negative values
	of $\widetilde D_{T,t}$ are significant). Then we may write
	\begin{equation} \label{DIST2}
	\widetilde D_{T,t}=\frac{1}{T^2} \sum_{j,k=1}^T \widetilde
	W_{j,k}+\frac{1}{(T+t)^2} \sum_{j,k=1}^{T+t} \widetilde
	W_{j,k}-\frac{2}{T(T+t)} \sum_{j=1}^T \sum_{k=1}^{T+t} \widetilde
	W_{j,k},
	\end{equation}
	where  $\widetilde W_{j,k}=\widetilde W(\Upsilon_j-\Upsilon_k)$ with
	\begin{equation} \label{szekely}
	\widetilde W(x)=\int_{\mathbb R^m} (1-\cos(u'x))w(u) du. \end{equation}
	If we let in (\ref{szekely}) $w(u)=\|u\|^{-(m+a)}, \ 0<a<2$, then \[
	\widetilde W(x)=\widetilde C \|x\|^{a},
	\] where $\widetilde C$ is a known constant depending only on $m$ and $a$,
	and hence in this case too our statistic comes in a closed-form expression
	suitable for computer implementation. Note that this weight function was first
	used by \citet{szekely05}, and later employed by \citet{matteson14}
	in change-point analysis.
	
	The choice for the weight function $w(\cdot)$ is usually based upon computational considerations. In fact if $w(\cdot)$ integrates to one (probably following some scaling) and satisfies $w(-u)=w(u)$ then the integral $W(\cdot)$ figuring in  (\ref{INT}) can be interpreted as the CF of a symmetric around zero random variable having density $w(\cdot)$. In this connection $w(\cdot)$ can be chosen as the density of any such distribution. Hence the choice $e^{-a\|u\|^2}$ corresponds to the multivariate normal density but for computational purposes any density with a simple CF will do. Other examples for $w(\cdot)$ are, for instance, the multivariate spherical stable density, the Laplace density, mixtures of normals or Laplace distributions, and combinations thereof. In fact, one might wonder whether there is a weight function which is optimal in some sense. The issue is still open but based on earlier finite-sample results it appears that the issue of the choice of $w$ is similar to the corresponding problem of choosing a kernel and a bandwidth in nonparametric estimation: Most weight functions (kernels) render similar behavior of the test statistic, but there is some sensitivity with respect to the ``bandwidth'' parameter $a>0$ figuring in (\ref{EXPW}). This is a highly technical problem that has been tackled only under the restrictive scenario of testing goodness-of-fit for a given parametric distribution, and even then a good choice of $a$ depends on the direction away from the null hypothesis; see \citet{tenreiro09}. Thus in our context the approach to the weight function is in some sense pragmatic: We use the Gaussian weight function which has become something of a standard, and investigate the behavior of the criterion over a grid of values of the weight parameter $a$. However, the alternative weight function below (\ref{szekely}) has also been tried, giving similar results.
	
	Another important user-specified parameter of our procedure is the order $m$ that determines the dimension of the joint distribution which is monitored for stationarity. Of course
	having a sample size $T$ of training data already imposes the obvious restriction $m\leq T$, but if $m$ is only slightly smaller than $T$, then we do not expect the ECF to be a reliable estimator of its population counterpart. The situation is similar to the problem of order-choice when estimating correlations from available data; see \citet{brockwell91}, p. 221, and \citet{box70}, p. 33, for some general guidelines. In our Monte Carlo results we only consider cases where $m\leq 4$ (very small compared to $T$), but it is reasonable to assume that we could let the dimension $m$ grow as $T$ increases; refer to Section 6.
	
	We close this section by noting that the statistic in (\ref{DIST}) compares the ECF computed from data up to time $T$,
	i.e.\ the ECF of the training data, with the ECF of all data available at the current monitoring time $T+t$, i.e.\ the data
	made available both before and after time $T$. Another option is to consider the statistic
	\begin{equation} \label{DIST3}
	D^*_{T,t}=\int_{\mathbb R^m}
	\left|\widehat \psi_T(u)-\widehat \psi^{(T+t)}(u)\right|^2 w(u) du,
	\end{equation}
	where
	\begin{equation} \label{ECF2}  \widehat \psi^{(J)}(u)=\frac{1}{J} \sum_{j=1}^J e^{i u' \Upsilon_{T+j}},
	\end{equation}
	in which case we compare the ECF of the training data to the ECF computed from the observations after time $T$, i.e.\ after
	monitoring has begun. Below we will often write $D_{T,t}$ as our statistic, but the methods also apply to the case of the other
	statistics considered in this section.  We finally close the section by noting that throughout we have assumed that $\{X_t\}$ is
	univariate. However extension to the multivariate context does not seem to present us with any undue complications.
	
	\section{Detector statistics and stopping rule}
	
	As already mentioned we consider on-line procedures whereby
	the test is applied sequentially on a dynamic data set which is steadily updated over time with the arrival of new observations. In this context, the null hypothesis is rejected when the value of a suitable detector statistic exceeds an appropriately chosen constant {\it{for the first
			time}}. Otherwise we continue monitoring.
	These statistics are commonly defined by a corresponding stopping
	rule. In order to define this stopping rule, and based on asymptotic
	considerations, we need to introduce an extra weight function in
	order to control the large-sample probability of type-I error. In
	particular we employ the detector statistic
	\begin{equation} \label{Detector}
	\Delta_{T,t}=\frac{1}{ q^2_\gamma\left(\frac{t}{T}\right)}
	\left(\frac{T+t-m+1}{\sqrt{T-m+1}}\right)^2 D_{T,t},
	\end{equation}
	where $D_{T,t}$ is defined by (\ref{DIST1}) and,
	\begin{equation} \label{weight}
	q_\gamma(s)=(1+s)\left(\frac{s}{s+1}\right)^\gamma, \ \gamma \in [0,1/2).
	\end{equation} Here $q_\gamma$ is an extra weight function needed to
	control (asymptotically) the probability of type-I error for the sequential test procedure. The parameter $\gamma$ figuring in (\ref{weight}) gives some flexibility to the resulting procedure. Specifically, if early violations are expected then the value of $\gamma$ should be close to 1/2, while values closer to zero are appropriate for detecting changes occurring at later stages; see \citet{aue06}.

	
	
	\medskip
	As already mentioned, it is clear that since the training data
	$\{X_1,...,X_T\}$ are assumed to involve no change, the monitoring
	period begins with time $t=T+1$. Typically this monitoring continues
	until time $T(L+1)$, where $L$ denotes a fixed integer, and if
	$L<\infty$ we call the corresponding procedure closed-end. Otherwise
	(i.e.\  if $L=\infty$), we have an open-end
	procedure. The corresponding stopping rule is specified as
	\begin{equation} \label{INAR12}
	\tau(T;L)= \tau(T)=
	\begin{cases}
	\inf\{1<t\leq L T: \Delta_{T,t} > c_\alpha\}, & \\
	+\infty, \textrm{ if }\Delta_{T,t} \leq c_\alpha \textrm{ for all } 1<t\leq L T,
	\end{cases}
	\end{equation}
	where $c_\alpha$ is a constant that guarantees that the test has size equal to $\alpha$, asymptotically.
	
	The main problem is then to find an approximation for the critical value $c_{\alpha}$ and to  investigate consistency of the test procedures. Particularly, we require that under
	${\cal{H}}_0$ and for a prechosen value of $0<\alpha<1$,
	\begin{equation} \label{level}
	\lim_{T\to\infty}   P_{{\cal{H}}_0}(\tau(T)<\infty)=\alpha,
	\end{equation}
	while under alternatives we want
	\begin{equation} \label{consist}
	\lim_{T\to\infty}   P(\tau(T)<\infty)=1.
	\end{equation}
	
	We close this section by noting that along the same lines one may
	suggest retrospective tests, i.e.\ tests in which all data have
	already arrived at the time of testing, and one is interested to see whether a structural
	change has occured over the time period leading to that time. To this end let $X_1,...,X_T$, denote not the training data
	(which involve no change), but data during the acquisition of which
	a structural change might have occurred. Then the appropriate
	statistic is that of (\ref{DIST3}) with the ECF $\widehat \psi_t(\cdot)$
	computed at each time $1\leq t\leq T$, with data $X_k, \
	1\leq k\leq t$, while the ECF  $\widehat \psi^{(t)}(\cdot) $ is computed
	in an analogous manner but with data  $X_k, \ t<k \leq T$.
	
	\section{Asymptotics}
	\smallskip
	In this section we study the limit behavior of the test procedure.
	Particularly, we have to study the limit behavior of
	\begin{equation} \label{MAX}
	\max_{1\leq t\leq LT}\Delta_{T,t}(\gamma),
	\end{equation}
	with $\Delta_{T,t}(\gamma)$ defined in (\ref{Detector}) under the
	null hypothesis  as well as under alternatives. The limit is always
	for $T\to\infty$ and $\tau$ fixed. The proofs of the lemmas and
	theorems in this Section are deferred to the Appendix.
	
	Define
	\begin{eqnarray}\label{PSS} S_{T, t} (u)= \frac{1}{\sqrt{T}}
	\sum_{j=1}^{T+t}f(\Upsilon_j, u)
	-\frac{1}{\sqrt{T}}\Big(\frac{T+t}{T}\Big)\sum_{j=1}^T f(\Upsilon_j
	, u)\end{eqnarray} and
	\begin{eqnarray*} S_{T, t}^m (u)=
		\frac{1}{\sqrt{T-m+1}} \sum_{j=m}^{T+t}f(\Upsilon_j,
		u)-\frac{1}{\sqrt{T-m+1}}\Big(\frac{T+t-m+1}{T-m+1}\Big)\sum_{j=m}^T
		f(\Upsilon_j , u),\end{eqnarray*} where $f$ is a two dimensional
	real-valued function on ${\mathbb R^m}\times{\mathbb R^m}$, defined
	by
	\begin{eqnarray}\label{function f}
	f(x , u )=(\cos( u^{'}x)-{\rm Re}\{\varphi_{\Upsilon_m}
	(u)\}, \sin (u^{'} x)-{\rm Im}\{ \varphi_{\Upsilon_m }(u)\})^{'},
	\end{eqnarray} we can express
	\begin{equation} \label{modified DIST}
	\Delta_{T,t}(\gamma)=\frac{ \int_{\mathbb R^m} ||S_{T, t}^m (u)||^2
		w(u) du}{ q^2_\gamma (\frac{t}{T})},
	\end{equation}
	where $S_{T,t}^m$ is replaced with $S_{T,t}$ when investigating the
	asymptotic behavior of $\Delta_{T,t}(\gamma)$.

	To obtain the limiting distribution of \eqref{MAX}, we consider
	this problem within a more general framework because the techniques
	used here would be applicable to other situations, for example, the
	monitoring procedure that uses the probability generating function
	approach as in \citet{hudecova15}. In what follows, $\{Y_t\}$ is
	assumed to be a $m (\geq 1)$-dimensional strongly mixing strictly
	stationary sequence with mixing order $\alpha(n)$, that is,
	$$\sup_{A\in {\cal F}_{-\infty}^k, B\in {\cal F}_{k+n}^{\infty}}|P(A\cap
	B)-P(A)P(B)|\leq \alpha(n)\searrow 0, $$ where ${\cal F}_{a}^b $
	denotes the sigma field generated by $Y_a, \ldots, Y_b$. Suppose
	that $f$ is any real-valued function on $\mathbb R^m\times \mathbb
	R^l$, $l\geq 1$. Later, a real vector function $f$ will be
	considered.
	
	In view of (\ref{PSS}), for any integer $T>m$,
	define
	\begin{eqnarray}\label{PSP} S_{T} (\theta, u)= \frac{1}{\sqrt{T}}
	\sum_{j=1}^{T+[\theta T]}f(Y_j, u)
	-\frac{1}{\sqrt{T}}\Big(\frac{T+[\theta T]}{T}\Big)\sum_{j=1}^T
	f(Y_j , u),\end{eqnarray} where $\theta$ is any number between 0 and
	$L$. 
	Without loss of generality, we assume $Ef(Y_1 , u)=0$
	because $f(Y_j , u)$ can be replaced with $f(Y_j , u
	)-Ef(Y_1, u )$ in $(\ref{PSP})$.
	
	Below, we impose the conditions:
	\begin{description}\label{conditions}
		\item{\bf (A1)} The $w(\cdot)$ is non-negative, continuous and
		bounded with $\int_{\mathbb R^l} w(u)du>0$.
		\item{\bf (A2)} $\alpha(n)\leq C n^{-\nu}$ for some
		$\nu>4$ for some $C>0$.
		\item{\bf (A3)} $||f||_\infty:=\sup_{x, u}|f(x, u )|
		<\infty$, and for some $K>0$, $
		E|f(Y_1, u)-f(Y_1, v)|^2\leq K
		||u-v||$.
		\item {\bf (A4)} $\sigma^2 (u):= \sum_{h=-\infty}^{\infty} Cov (f(Y_1, u), f(Y_{1+h} , u
		))$ is continuous in $u$, and  $0<\inf_{u}
		\sigma(u)\leq \sup_{u}\sigma(u)<\infty$.
	\end{description}
	
	Condition {\bf (A2)} is quite mild since a broad class of time
	series sequences, including ARMA and GARCH processes, are strongly
	mixing with order decaying to 0 geometrically. Owing to Lemma 2.1 of \citet{kuelbs80}, we have that for all $\zeta>1$,
	\begin{eqnarray}\label{covariance inequality}
	|Cov (f(Y_1 , u ), f(Y_{1+h}, u )|\leq
	12\alpha(h)^{1/\zeta}| ||f||_\infty^2,
	\end{eqnarray}
	so that $\sigma^2 (u)$ in {\bf (A4)} converges absolutely.
	Moreover, due to Theorem A and Theorem 2.1 of \citet{Yang2007},
	Conditions {\bf (A1)}-{\bf (A4)} particularly guarantee the maximal moment inequality for
	$S_k=\sum_{j=1}^k f(Y_j, u)$, namely, for any $r\geq 2$, there exists $K_r>0$ depending only on $r$ and  $||f||_\infty$
	such that
	\begin{eqnarray}\label{maximal inequality}
	\sup_u E\max_{1\leq k\leq T}|S_k |^r\leq K_r T^{r/2}\ {\rm for \ all }\ T\geq 1.
	\end{eqnarray}
	Then, we have the following:
	\medskip\\ \indent
	{\bf Lemma 1}. {\it Under {\bf (A1)}-{\bf (A3)},  as
		$T\rightarrow\infty$,
		\begin{eqnarray*}
			Z_{T}(\theta, u):=\frac{S_{T} (\theta, u)}{
				q\Big(\frac{[T\theta]}{T}\Big)}\stackrel{d}{\longrightarrow}Z(\theta,
			u):=\sigma(u){\cal B} \Big(\frac{\theta}{1+\theta}\Big)/
			\Big(\frac{\theta}{1+\theta}\Big)^\gamma,
		\end{eqnarray*}
		where $Z_{T}(\theta, u)=0$ when $\theta<1/T$ and ${\cal B}$ denotes
		a standard Brownian motion. Further, $$(Z_{T}(\theta_1, u_1
		),\ldots, Z_{T}(\theta_k, u_k ))\stackrel{d}{\longrightarrow}
		(Z(\theta_1, u_1 ),\ldots, Z(\theta_k, u_k))$$ for any $\theta_i\in
		[0, L]$ and $u_i\in {\mathbb R^l}, i=1,\ldots, k$, $k\geq 1$.}
	\medskip\\ \indent
	Based on this, we have the following weak convergence results \citep[cf.][]{billingsley68}:
	\medskip\\ \indent
	{\bf Theorem 1}. {\it
		Let
		$$ {\cal Z}_T (\theta)=\int_{\mathbb R^l} Z_{T}^2 (\theta, u) w(u)du \ \ {\rm
			and}\ \ {\cal Z}(\theta)=\frac{|{\cal
				B}(\frac{\theta}{1+\theta})|^2}{(\frac{\theta}{1+\theta})^{2\gamma}}\int_{\mathbb R^l}
		\sigma^2 (u) w(u)du.$$
		Assume that $0\leq \eta\leq L$ is any real number that can take the value of 0 only when $\gamma=0$.  Then, under {\bf (A1)}-{\bf (A4)},
		as $T\rightarrow\infty$, $$\sup_{\eta \leq \theta\leq L}{\cal Z}_T
		(\theta) \stackrel{d}{\rightarrow} \sup_{\eta\leq\theta\leq L}{\cal
			Z}(\theta),$$ namely,\
		\begin{eqnarray}\label{Theorem1c}
		\sup_{\eta\leq \theta\leq L}\int_{\mathbb R^l} \frac{S_{T}^2 (\theta, u)}{
			q^2_\gamma\Big(\frac{[T\theta]}{T}\Big)} w(u)du
		\stackrel{d}{\longrightarrow}\sup_{\frac{\eta}{1+\eta}\leq s\leq
			\frac{L}{1+L}}\frac{|{\cal B}(s)|^2}{s^{2\gamma}}\int_{\mathbb R^l}
		\sigma^2(u)w(u)du,
		\end{eqnarray}
		and thus,
		\begin{eqnarray}\label{Theorem1d}
		\max_{\eta T \leq k \leq LT}\int_{\mathbb R^l}\frac{S_{T}^2 (\frac{k}{T}, u)}{
			q^2_\gamma(\frac{k}{T})} w(u)du
		\stackrel{d}{\longrightarrow}\sup_{\frac{\eta}{1+\eta}\leq s\leq
			\frac{L}{1+L}}\frac{|{\cal B}(s)|^2}{s^{2\gamma}}\int_{\mathbb R^l}\sigma^2(u)
		w(u)du.
		\end{eqnarray}
	}
	
	{\bf Remark 1}. {\it
		(a) The result shows that  when $\gamma=0$, $\sup_{0\leq \theta\leq L}{\cal Z}_T
		(\theta) \stackrel{d}{\rightarrow} \sup_{0\leq\theta\leq L}{\cal
			Z}(\theta)$. However, the range of $\theta$ for $\gamma>0$ is not allowed to be $[0, L]$ but $[\eta, L]$ for  $0<\eta<L$. Although the result on $[0, L]$ can be easily claimed, it is not easy to verify. The difficulty lies in checking the tightness (equicontinuity) condition such as that in Theorem 15.6 of \citet{billingsley68}.
		The approaches taken by \citet{aue06} and \citet{hudecova15,hudecova15b} are not directly applicable to solving this problem.
		In practice, though, the $\eta$ can be chosen arbitrarily small, which can be justified from a practical viewpoint that
		a change is assumed to occur at $t=T+[T\tau]$ for some $0<\tau<1$.

		(b) The $\sigma^2 (u)$ can be estimated by
		\begin{eqnarray}\label{variance estimate}  \hat{\sigma}_T^2 (u)= \sum_{|h|\leq h_T}
		\frac{1}{T-|h|}\sum_{j=1}^{T-|h|} \{f(Y_1, u)-\hat\mu_T
		(u)\}\{f(Y_{1+h}, u )-\hat\mu_T (u)\},\end{eqnarray} where
		$\hat\mu_T=\frac{1}{T}\sum_{j=1}^T f(Y_j , u )$ and $\{h_T\}$ is a
		sequence such that $h_T / T \rightarrow 0$ as $T\rightarrow\infty$,
		for example, $h_T = T^{1/3}$.
		
		(c) If {\bf (A3)} and {\bf (A4)} hold with $u, v\in \mathbb R^l$ replaced by $u, v$ in a compact subset  $U$  of ${\mathbb R^l}$,
		then the results in (\ref{Theorem1c}) and
		(\ref{Theorem1d}) hold when the integral over ${\mathbb R^l}$ is replaced with
		the integral over $U$.
	}
	
	{\bf Remark 2}. {\it In our study the strong mixing condition is assumed and the strong approximation theorems are used to achieve the C.L.T. of the relevant partial sum process. However, one can view $S_T (\theta, u)$ as an element of $D[0,1+L]$ by resetting $\vartheta=1+\theta$ and using
		\begin{eqnarray}\label{tildePSP} \tilde S_{T} (\vartheta, u)= \frac{1}{\sqrt{T}}
		\sum_{j=1}^{[T\vartheta]}f(Y_j, u)
		-\frac{1}{\sqrt T}\frac{[T\vartheta]}{T}\sum_{j=1}^T
		f(Y_j , u).\end{eqnarray}
		Our test statistic can be redefined with $\{\tilde S_T (\vartheta, u); 1\leq \vartheta\leq 1+L\}$, which is a subprocess of $\{\tilde S_T (\vartheta, u); 0\leq \theta\leq 1+L\}$, and thereby, the strong approximation theorems may not be required. 
		In our simulation study, the stationary bootstrap method is used to calculate critical values.
		A proof for the consistency of the bootstrap method is provided in Supplementary material.
		
	}
	

	The above result is useful to test the hypotheses:
	\medskip\\ \indent
	${\cal K}_0$:  No change of $E f(Y_k, u)$ over
	$k=T+1,\ldots, T+LT$\ \ vs.\ \ ${\cal K}_1$: not ${\cal K}_0$.
	\medskip\\ \noindent In
	order to conduct a test, one can employ the test: for some nonnegative integer $\eta$ (=0 when $\gamma=0$ and $>0$ when $\gamma>0$),
	$$ {\cal T}(T) =\max_{\eta T\leq k \leq LT}\int_{\mathbb R^l} \frac{S_{T}^2 (\frac{k}{T},
		u)}{ q^2_\gamma(\frac{k}{T})} w(u)du.$$  Based on Theorem 1, we
	reject ${\cal K}_0$ at the significance level $\alpha$ if ${\cal
		T}(T)\geq c_\alpha$, where $c_\alpha$ is the number such that
	$$P\left(\sup_{\frac{\eta}{1+\eta}\leq s\leq
		\frac{L}{1+L}}\frac{|{\cal B}(s)|^2}{s^{2\gamma}}\int_{\mathbb R^l}
	\sigma^2 (u)w(u)du\geq c_\alpha\right) =\alpha.$$ In practice, the
	$c_\alpha$ can be obtained through Monte Carlo simulations using an
	estimate $\hat\sigma_T^2(u)$ of $\sigma^2(u)$.
	\medskip\\ \indent
	Next, we extend the above results to  $d$-dimensional functions
	$f=(f_1, \ldots, f_d )$, $d\geq 1$. Let $\Sigma(u)$ denote  the
	positive definite $d\times d$ matrix whose $(i,j)$th entry
	$\Sigma_{ij}(u)$ is $\sum_{h=-\infty}^\infty Cov (f_i (Y_1, u), f_j (Y_{1+h}, u ))$. Below, we
	assume that
	\begin{description}
		\item {\bf (A3)}$^{'}$ All $f_i$ enjoy the properties in {\bf (A3)}.
		\item{\bf (A4)}$^{'}$ $u\rightarrow \Sigma(u)$ is continuous and $0<\inf_u |\Sigma (u)||\leq
		\sup_{u}||\Sigma (u)||<\infty$.
	\end{description}
	Then, we have the following:
	\medskip\\ \indent
	{\bf Theorem 2}. {\it Suppose that {\bf (A1)}, {\bf (A2)}, {\bf
			(A3)}$^{'}$ and {\bf (A4)}$^{'}$ hold. Let $\eta$ be the one in Theorem 1. Then, as $T\rightarrow\infty$,
		\begin{eqnarray}\label{weak conv2}
		\sup_{\eta\leq \theta\leq L}\int_{\mathbb R^l}
		\frac{||S_{T}(\theta,u)||^2}{q^2_\gamma\Big(\frac{[T\theta]}{T}\Big)}
		w(u)du \stackrel{d}{\longrightarrow}\sup_{\frac{\eta}{1+\eta}\leq s\leq
			\frac{L}{1+L}}\frac{1}{s^{2\gamma}}{\cal B}_d^{'}(s)\Big(\int_{\mathbb R^l}
		\Sigma (u)w(u)du\Big ) {\cal B}_d (s),
		\end{eqnarray}
		where ${\cal B}_d$ is a standard $d$-dimensional Brownian motion,
		and thus,
		\begin{eqnarray}\label{weak conv22}
		\sup_{\eta T\leq k\leq LT}\int_{\mathbb R^l}
		\frac{||S_{T}(\frac{k}{T},u)||^2}{q^2_\gamma\Big(\frac{k}{T}\Big)}
		w(u)du \stackrel{d}{\longrightarrow}\sup_{\frac{\eta}{1+\eta}\leq s\leq
			\frac{L}{1+L}}\frac{1}{s^{2\gamma}}{\cal B}^{'}_d(s)\Big(\int_{\mathbb R^l}
		\Sigma (u)w(u)du\Big ) {\cal B}_d (s).
		\end{eqnarray}   }
	%
	
	The testing procedure based on the above result is similar to the
	univariate case. That is, we reject the null of no change if
	\begin{eqnarray}\label{test2}
	{\cal T}(T)=\sup_{\eta T\leq k \leq LT}\int_{\mathbb R^l}
	\frac{||S_{T}(\frac{k}{T},u)||^2}{q^2_\gamma\Big(\frac{k}{T}\Big)}
	w(u)du\geq c_\alpha,
	\end{eqnarray}
	where $c_\alpha$ satisfies
	\begin{eqnarray}\label{test2 critical}
	P\left(\sup_{\frac{\eta}{1+\eta}\leq s\leq \frac{L}{1+L}}\frac{1}{s^{2\gamma}}{\cal
		B}^{'}_d (s)\Big(\int_{\mathbb R^l} \Sigma (u)w(u)du\Big ) {\cal
		B}_d (s)\geq c_\alpha\right )=\alpha.
	\end{eqnarray}
	The $c_\alpha$ can be obtained through  Monte Carlo simulations
	using a suitable estimate $\hat\Sigma_T (u)$ of $\Sigma(u)$.
	\medskip\\ \indent
	Now, we are ready to apply Theorem 2 with $d=2$ and $l=m$ to $\max_{\eta T\leq
		t\leq LT}\Delta_{T,t}(\gamma)$ in $(\ref{MAX})$. Assume that
	$\{X_t\}_{t\in \mathbb N}$ is a strongly mixing process with mixing
	order $\alpha(n)$ satisfying ${\bf (A2)}$ and $w(\cdot)$ satisfies
	{\bf (A1)}. Note that under ${\cal H}_0$ in $(\ref{NULL})$,
	\begin{eqnarray}\label{stopping rule}
	\max_{\eta T\leq t\leq
		LT}\Delta_{T,t}(\gamma)\stackrel{d}{\longrightarrow}\sup_{\frac{\eta}{1+\eta}\leq
		s\leq \frac{L}{1+L}}\frac{1}{s^{2\gamma}}{\cal B}^{'}_2
	(s)\Big(\int_{\mathbb R^m} \Sigma (u)w(u)du\Big ) {\cal B}_2 (s),
	\end{eqnarray} where $\Sigma(u)$ is a $2\times 2$ matrix whose
	$(ij)$th entries are $\Sigma_{11}(u)= Var (\cos( u^{'}
	\Upsilon_{1}))$, $\Sigma_{22}(u)= Var ( \sin( u^{'} \Upsilon_{1}))$,
	and $\Sigma_{12}(u)=\Sigma_{21}(u)=Cov (\cos( u^{'} \Upsilon_{1}),
	\sin (u^{'} \Upsilon_{1}))$ (see $(\ref{function f})$) and satisfies
	$\sup_{u\in {\mathbb R^m} } ||\Sigma (u)||<\infty$. We reject ${\cal
		H}_0$ in ($\ref{NULL})$, if $\max_{1\leq t\leq
		LT}\Delta_{T,t}(\gamma)\geq c_\alpha$ where $c_\alpha$ is the number
	in $(\ref{test2 critical})$.
	
	Meanwhile, to examine the consistency of the monitoring procedure,
	we consider the alternative hypothesis: $$\tilde{\cal H}_1 :
	\psi_{\Upsilon_1}(u_1 )=\cdots=\psi_{\Upsilon_{[\tau
			T]}}(u_1)=\psi_1 \neq\psi_2= \psi_{\Upsilon_{[\tau T]+1}}(u_1
	)=\cdots= \psi_{\Upsilon_{LT}}(u_1)$$ for some $1<\tau<L$ and
	$u_1\in{\mathbb R^m}$. Then, in view of $(\ref{PSS})$, we have that
	under $\tilde{\cal H}_1$,
	\begin{eqnarray*}\label{PS} \frac{S_{T, [\tau T]}
			(u_1)}{T}\rightarrow \tau (\psi_1-\psi_2)\end{eqnarray*} in
	probability, so that $ \max_{1\leq t\leq LT}\Delta_{T,[\tau
		T]}(\gamma)\rightarrow\infty$ in probability, which asserts the
	consistency of our monitoring procedure as in $(\ref{level})$.
	
	The matrix $\Sigma(u)$ should be estimated properly as done in
	($\ref{variance estimate})$. Also, as mentioned earlier in $(\ref{DIST1})$
	and $(\ref{EXPW})$, $\Delta_{T, t}$ can be expressed in closed form
	for particular weight functions.
	By way of example in the empirical study below, we
	will illustrate the performance of the monitoring procedure with
	the Gaussian weight function.
	
	\section{The resampling procedure}\label{sec:resampling.procedure}
	As already seen in the previous section, the asymptotic distribution of the detector statistic
	in \eqref{MAX} under the null hypothesis ${\cal{H}}_0$ depends on factors that are unknown in our entirely nonparametric context.
	For this reason we employ the stationary bootstrap procedure
	\citep[see][]{politis1994} in order to estimate
	the critical value $c_\alpha$ of the test. To be concrete, the estimator $\hat c_\alpha$ for $c_\alpha$ is given by the relation
	\[
	\Prob\Big(\max_{1\le t\le LT} \Delta^*_{T,t} \ge \hat c_\alpha \,\Big|\, X_1,\ldots,X_T \Big) = \alpha,
	\]
	where $\Delta^*_{T,t} = \Delta_{T,t}(X_1^*,\ldots,X_{T+t}^*)$ denotes the statistic in \eqref{Detector} applied to the resampled observations $X_1^*,\ldots,X_{T+t}^*$ obtained using the resampling scheme below.
	Notice that $c_\alpha$ is estimated by making use of only the training data, i.e.\  all data available at time $T$.
	
	First, wrap the observations around a circle so that for $k>T$, $X_k$ is defined to be $X_{k\bmod{T}}$ with $X_0 = X_T$.
	Given a block size $\ell$, define the overlapping blocks $B_{k,\ell} = \{X_k,\ldots,X_{k+\ell-1}\}$, $k = 1,\ldots,T$.
	Let $N=T(1+L)$ and proceed as follows:
	\begin{enumerate}
		\item Independently of the $X_k$, generate independent observations $\ell_1,\ell_2,\ldots$ from the geometric distribution with parameter $p_B$. Let $n_\ell=\min\{n:\sum_{k=1}^n \ell_k \ge N\}$.
		\item Independently of the $X_k$ and the $\ell_k$, generate independent observations $I_1,I_2,\ldots,I_{n_\ell}$ from the discrete uniform distribution on $\{1,\ldots,T\}$. Define the sequence of bootstrap blocks $\{B_{k,\ell_k}^*\}_{k=1}^{n_\ell}$ by $B_{k,\ell_k}^* = B_{I_k,\ell_k}$.
		\item Concatenate the $B_{k,\ell_k}^*$ and select the first $N$ observations as the bootstrap sample $X_1^*,\ldots,X_N^*$.
		\item Treat the first $T$ bootstrap observations as a training sample
		and calculate $\Delta^*_{T,t} = \Delta_{T,t}(X_1^*,\ldots,X_{T+t}^*)$ for each
		$t = 1,\ldots,LT$.
		\item Calculate $M^* = \max_{1\le t \le LT} \Delta^*_{T,t}$.
	\end{enumerate}
	Repeat steps 1--5 a large number of times, say $B$, to obtain the ordered statistics $M_{(1)}^* \le \cdots \le M_{(B)}^*$.
	An approximate value for $c_\alpha$ is then given by $M_{(\lfloor B(1-\alpha) \rfloor)}^*$,
	where $\lfloor x \rfloor$ denotes the floor of $x\in\mathbb{R}$.
	
	In order to choose an appropriate value for $p_B$ appearing in step~1 above, we employ the selection procedure proposed by \citet{politis2004}; also see \citet{patton2009}. This choice of $p_B$ is used throughout the simulations discussed in the following section.
	
	\section{Monte Carlo results}
	We investigate the performance of the monitoring procedure defined by (\ref{Detector})
	and (\ref{MAX}) with criterion given by (\ref{DIST1}),
	and weight function $w(u)=e^{-a\|u\|^2}, \ a>0$.
	The results corresponding to criterion (\ref{DIST2})
	are similar and therefore are not reported.
	We consider  the following data generating processes (DGPs):
	\begin{align*}
	&\text{DGP.S1: } X_t = \varepsilon_t, \{\varepsilon_{t}\} \stackrel{i.i.d.}{\sim} \text{N(0,1)} \\
	&\text{DGP.S2: } X_t = 0.5X_{t-1} + \varepsilon_t, \\
	&\text{DGP.S3: } X_t = h_t \varepsilon_t, \text{ with } h_t^2 = 0.2 + 0.3X_{t-1}^2, \\
	&\text{DGP.S4: } X_t = h_t \varepsilon_t, \text{ with } h_t^2 = 0.1 + 0.3X_{t-1}^2 + 0.3h_{t-1}^2, \\
	&\text{DGP.S5: } X_t = h_t \varepsilon_t, \text{ with } h_t^2 = 0.1 + 0.7X_{t-1}^2 + 0.3h_{t-1}^2, \\
	&\text{DGP.S6: } X_t = \beta_t X_{t-1} + \varepsilon_t, \beta_t = 0.5\beta_{t-1} + \eta_t, \{\eta_{t}\} \stackrel{i.i.d.}{\sim} N(0,0.1^2), \\
	&\text{DGP.S7: } X_t = \eta_t, \{\eta_{t}\} \stackrel{i.i.d.}{\sim} \text{Cauchy}, \\
	&\text{DGP.P1: } X_t = \varepsilon_t + \I\{t > V\}, \\
	&\text{DGP.P2: } X_t = \varepsilon_t (1 + \I\{t > V\}), \\
	&\text{DGP.P3: } X_t = (1+\sqrt{2}\varepsilon_t) \I\{t \le V\} + \varepsilon_t^2 \I\{t > V\}, \\
	&\text{DGP.P4: } X_t = \varepsilon_t \I\{t \le T\} + \varepsilon_t \exp(\tfrac{1}{2}-|\tfrac{1}{2} - (T-t)/LT|) \I\{t > T\}, \\
	&\text{DGP.P5: } X_t = \varepsilon_t \I\{t \le V\} + \xi_t \I\{t > V\}
	\end{align*}
	where $X_0 = 0$, and $V=T(1+UL)$, with $U$ uniformly distributed on $(0,\tfrac{4}{5})$.
	The sequence $\{\xi_t\}_{t \in \mathbb{N}}$ is an independent sequence of iid $\alpha$-stable random variables
	with characteristic function $\exp\{-|t|(1+0.5i\sgn(t)\log|t|/\pi)\}$.
	
	The DGPs consist of a range of serial dependence structures, such as AR, ARCH and GARCH structures.
	The first seven processes, DGP.S1--DGP.S7, satisfy the null hypothesis of strict stationarity
	and are introduced to study the size of our monitoring procedure.
	DGP.S1 and DGP.S7 consists of i.i.d.\ observations, whereas DGP.S2--DGP.S6 introduce some dependence structure without violating the null hypothesis.
	To examine the power of the procedure in different settings, DGP.P1--DGP.P5
	remain stationary throughout the training period, after which a break in stationarity
	is introduced during monitoring.
	DGP.P1 and DGP.P2 introduce breaks in the first and second moments (respectively) of the process,
	whereas DGP.P3 introduces a change in the distribution without affecting the first two
	moments of the DGP (i.e.\ the process remains weakly stationary throughout the monitoring period).
	
	We compare the performance of our test against that of several related tests proposed recently in the literature. We consider the following tests:
	\begin{itemize}
		\item a kernel-based nonparametric test for changes in the marginal density function proposed by \citet{lee2004}, denoted by $\Delta^{(L)}$;
		\item a related kernel-based test of \citet{kapetanios09}, denoted by $\Delta^{(K)}$;
		\item a Cram\'er--von~Mises type test for changes in distribution based on the empirical distribution function, proposed by \citet{ross2012}, denoted by $\Delta^{(R)}$.
	\end{itemize}
	Although the tests $\Delta^{(L)}$ and $\Delta^{(K)}$ were originally proposed as offline (retrospective) tests for stationarity, they have been adapted to online monitoring procedures.
	A similar adaptation of the retrospective test by \citet{hong2016} to an online procedure was also considered, but the obtained simulation results (not shown) suggest that the adapted procedure does not respect the desired nominal size. We therefore exclude this procedure from our study.
	
	The values in Tables~\ref{tab:S1}--\ref{tab:Pother} represent the percentage of times
	that $\mathcal{H}_0$ was rejected out of 1\,000 independent Monte Carlo repetitions
	for the various DGPs. For our test, denoted by $\Delta_a$ in the tables,
	we report the results for various choices of the parameters $m$ and $a$ appearing in the test.
	To estimate the critical value of each test we employed the bootstrap scheme
	described in Section~\ref{sec:resampling.procedure} and, since the calculations are very
	time consuming, made use of the warp-speed method of \citet{giacomini13}.
	Throughout, we took $\alpha=5\%$.
	
	An overall assessment of the figures in Tables~\ref{tab:S1}--\ref{tab:Pother}
	brings forward certain desirable features of the method:
	a reasonable degree of approximation of the nominal level,
	as well as satisfactory power against alternatives.
	Although our test exhibits some degree of moderate size distortion with the
	autoregressive process DGP.S2, this over-rejection certainly diminishes
	with increasing sample sizes. Note however that overall the nominal level is respected
	reasonably well, and also in some cases where some of the other considered tests are substantially over-sized.
	
	As expected, the power of the tests increases with the sample size,
	and also with the extent of violation of the
	null hypothesis of strict stationarity.
	Overall, in terms of power, our test performs reasonably
	well when compared to the other existing procedures and even outperforms
	these tests for many of the DGPs considered.
	A further factor influencing the percentage of rejection is the value of
	the weight parameter $a$, and in this respect it
	seems that an intermediate value $0.5\leq a\leq 1.5$ is preferable
	in terms of level accuracy as well as power,
	at least in most cases. 
	
	\begin{sidewaystable}
		\caption{Size results for DGP.S1--DGP.S4 for the procedure based on $\Delta_a$.}\label{tab:S1}
		\centering\scriptsize
		\addtolength{\tabcolsep}{-2pt}
		\begin{tabular}{cccccccccccccccccccccccccccccccccccccc}
			\toprule
			&&& \multicolumn{5}{c}{DGP.S1} & \multicolumn{5}{c}{DGP.S2} & \multicolumn{5}{c}{DGP.S3} & \multicolumn{5}{c}{DGP.S4} \\
			\cmidrule(lr){4-8} \cmidrule(lr){9-13} \cmidrule(lr){14-18} \cmidrule(lr){19-23}
			$m$ & $T$ & $L$
			& \multicolumn{1}{c}{$\Delta_{0.1}$} & \multicolumn{1}{c}{$\Delta_{0.5}$} & \multicolumn{1}{c}{$\Delta_{1}$} & \multicolumn{1}{c}{$\Delta_{1.5}$} & \multicolumn{1}{c}{$\Delta_{5}$}
			& \multicolumn{1}{c}{$\Delta_{0.1}$} & \multicolumn{1}{c}{$\Delta_{0.5}$} & \multicolumn{1}{c}{$\Delta_{1}$} & \multicolumn{1}{c}{$\Delta_{1.5}$} & \multicolumn{1}{c}{$\Delta_{5}$}
			& \multicolumn{1}{c}{$\Delta_{0.1}$} & \multicolumn{1}{c}{$\Delta_{0.5}$} & \multicolumn{1}{c}{$\Delta_{1}$} & \multicolumn{1}{c}{$\Delta_{1.5}$} & \multicolumn{1}{c}{$\Delta_{5}$}
			& \multicolumn{1}{c}{$\Delta_{0.1}$} & \multicolumn{1}{c}{$\Delta_{0.5}$} & \multicolumn{1}{c}{$\Delta_{1}$} & \multicolumn{1}{c}{$\Delta_{1.5}$} & \multicolumn{1}{c}{$\Delta_{5}$} \\
			\midrule
			1 & 100 & 1 & 3.6 & 3.8 & 3.8 & 4.0 & 4.1 & 8.4 & 9.6 & 9.6 & 9.7 & 9.9 & 7.8 & 5.3 & 6.0 & 5.4 & 4.3 & 3.8 & 4.3 & 3.4 & 2.9 & 3.6 \\
			& & 2 & 4.2 & 4.6 & 5.0 & 4.8 & 5.2 & 8.7 & 8.1 & 9.1 & 9.2 & 9.0 & 7.1 & 7.1 & 6.1 & 5.2 & 4.8 & 4.5 & 6.4 & 6.3 & 6.1 & 5.7 \\
			& & 3 & 6.1 & 4.1 & 4.5 & 4.2 & 3.9 & 7.4 & 8.5 & 8.4 & 8.4 & 8.4 & 7.9 & 7.8 & 6.6 & 6.4 & 6.6 & 8.9 & 9.5 & 7.6 & 7.1 & 7.0 \\ 
			\cmidrule{2-23}
			& 200 & 1 & 4.1 & 4.4 & 5.0 & 4.5 & 4.0 & 6.5 & 7.0 & 7.2 & 7.4 & 7.7 & 4.8 & 5.5 & 5.6 & 5.5 & 4.8 & 4.4 & 4.5 & 5.2 & 5.0 & 5.0 \\
			& & 2 & 5.1 & 4.3 & 4.2 & 4.2 & 4.4 & 7.4 & 7.7 & 7.7 & 7.5 & 7.4 & 5.7 & 6.0 & 5.5 & 5.2 & 4.9 & 6.9 & 7.2 & 8.5 & 7.7 & 6.0 \\
			& & 3 & 3.5 & 4.2 & 3.8 & 4.8 & 5.6 & 8.4 & 8.6 & 7.9 & 7.8 & 7.6 & 7.3 & 7.8 & 7.3 & 7.5 & 7.0 & 5.6 & 5.0 & 5.3 & 5.4 & 5.9 \\ 
			\cmidrule{2-23}
			& 300 & 1 & 4.5 & 4.7 & 4.1 & 3.6 & 3.7 & 6.6 & 7.9 & 7.9 & 7.8 & 7.9 & 5.7 & 6.0 & 6.4 & 5.1 & 3.7 & 5.7 & 5.2 & 5.2 & 5.1 & 5.1 \\
			& & 2 & 4.0 & 2.8 & 2.3 & 2.7 & 2.7 & 6.5 & 7.7 & 7.5 & 7.5 & 7.5 & 5.2 & 6.5 & 6.2 & 6.3 & 5.5 & 8.2 & 7.0 & 6.5 & 5.8 & 5.9 \\
			& & 3 & 7.3 & 7.1 & 7.7 & 7.4 & 6.9 & 7.5 & 7.7 & 7.6 & 7.9 & 7.7 & 4.4 & 4.6 & 5.8 & 6.2 & 5.2 & 6.7 & 6.9 & 6.4 & 6.2 & 4.7 \\ 
			\midrule
			2 & 100 & 1 & 3.7 & 4.4 & 3.8 & 3.6 & 4.3 & 7.5 & 9.1 & 9.1 & 9.1 & 9.6 & 8.0 & 6.3 & 5.5 & 5.1 & 4.0 & 3.4 & 3.7 & 4.2 & 3.2 & 3.6 \\
			& & 2 & 4.2 & 5.0 & 4.8 & 4.4 & 4.5 & 6.9 & 9.2 & 9.8 & 9.7 & 9.3 & 6.3 & 5.4 & 6.0 & 5.6 & 4.7 & 5.0 & 6.0 & 6.1 & 6.0 & 4.2 \\
			& & 3 & 3.6 & 4.4 & 4.6 & 4.4 & 4.4 & 6.7 & 8.4 & 8.5 & 8.2 & 8.5 & 5.9 & 6.5 & 6.9 & 6.7 & 5.9 & 7.0 & 7.0 & 6.7 & 6.5 & 7.2 \\ 
			\cmidrule{2-23}
			& 200 & 1 & 3.8 & 4.8 & 4.4 & 4.3 & 4.1 & 6.6 & 7.9 & 8.4 & 8.5 & 9.2 & 4.5 & 4.3 & 4.2 & 4.8 & 4.6 & 4.8 & 4.1 & 3.7 & 4.7 & 4.3 \\
			& & 2 & 5.0 & 3.7 & 4.3 & 4.1 & 4.3 & 6.0 & 7.0 & 7.5 & 7.7 & 7.9 & 4.8 & 5.3 & 5.3 & 5.1 & 4.8 & 5.9 & 6.9 & 7.1 & 6.9 & 5.7 \\
			& & 3 & 3.3 & 4.8 & 4.2 & 4.6 & 5.8 & 7.1 & 7.8 & 8.3 & 8.1 & 7.6 & 6.2 & 6.5 & 7.4 & 7.4 & 6.5 & 5.0 & 4.6 & 4.5 & 4.7 & 5.7 \\ 
			\cmidrule{2-23}
			& 300 & 1 & 4.0 & 3.7 & 4.7 & 4.4 & 3.6 & 5.7 & 7.5 & 7.6 & 7.9 & 7.8 & 5.4 & 4.9 & 4.7 & 4.5 & 3.3 & 6.1 & 5.5 & 5.5 & 5.5 & 5.0 \\
			& & 2 & 4.2 & 3.9 & 2.1 & 2.3 & 2.3 & 5.6 & 7.3 & 7.5 & 7.4 & 7.6 & 5.2 & 5.7 & 5.3 & 5.3 & 5.1 & 6.7 & 6.7 & 6.0 & 5.9 & 5.5 \\
			& & 3 & 6.6 & 7.3 & 6.4 & 7.1 & 7.1 & 6.8 & 7.7 & 7.3 & 7.1 & 7.6 & 4.7 & 4.1 & 4.2 & 5.1 & 5.3 & 7.0 & 5.9 & 5.6 & 6.0 & 4.5 \\ 
			\midrule
			4 & 100 & 1 & 2.8 & 4.9 & 4.9 & 5.1 & 4.9 & 2.7 & 8.1 & 8.6 & 8.8 & 9.2 & 5.3 & 6.0 & 5.8 & 5.9 & 4.6 & 3.5 & 4.4 & 4.0 & 4.1 & 3.3 \\
			& & 2 & 3.3 & 5.2 & 5.2 & 5.1 & 4.7 & 2.7 & 8.0 & 8.9 & 9.2 & 9.6 & 4.9 & 5.7 & 5.5 & 5.9 & 5.6 & 4.9 & 6.2 & 6.3 & 6.1 & 4.3 \\
			& & 3 & 2.8 & 4.3 & 4.7 & 5.0 & 5.2 & 2.6 & 7.5 & 8.2 & 8.4 & 8.2 & 5.3 & 6.3 & 5.7 & 5.9 & 5.0 & 6.5 & 6.6 & 6.4 & 6.2 & 6.1 \\ 
			\cmidrule{2-23}
			& 200 & 1 & 4.2 & 4.5 & 4.3 & 4.7 & 5.1 & 2.5 & 7.1 & 7.6 & 8.1 & 8.6 & 4.7 & 5.5 & 5.5 & 5.0 & 4.4 & 6.1 & 5.3 & 4.8 & 4.3 & 3.6 \\
			& & 2 & 3.6 & 4.4 & 4.9 & 5.2 & 5.5 & 2.6 & 6.4 & 6.8 & 7.4 & 7.7 & 5.3 & 5.4 & 5.2 & 4.8 & 4.2 & 5.6 & 7.0 & 7.4 & 7.1 & 5.2 \\
			& & 3 & 3.6 & 4.6 & 4.7 & 4.8 & 5.0 & 3.2 & 6.9 & 7.5 & 8.0 & 7.7 & 5.3 & 6.3 & 5.8 & 5.6 & 5.1 & 6.7 & 5.6 & 4.8 & 4.3 & 5.4 \\ 
			\cmidrule{2-23}
			& 300 & 1 & 3.7 & 4.8 & 4.9 & 4.9 & 4.9 & 2.6 & 6.3 & 6.9 & 7.5 & 7.7 & 6.1 & 6.2 & 5.8 & 5.6 & 4.4 & 6.2 & 7.4 & 6.0 & 5.4 & 5.2 \\
			& & 2 & 4.6 & 5.1 & 5.2 & 4.8 & 4.7 & 2.9 & 6.8 & 7.4 & 7.4 & 7.5 & 5.2 & 5.5 & 5.8 & 5.6 & 4.7 & 7.5 & 6.5 & 6.8 & 7.0 & 5.8 \\
			& & 3 & 5.1 & 5.4 & 5.3 & 5.5 & 5.5 & 2.7 & 6.9 & 7.6 & 7.3 & 7.2 & 5.7 & 5.9 & 5.7 & 5.7 & 5.1 & 5.2 & 4.9 & 4.7 & 5.2 & 5.7 \\
			\bottomrule
		\end{tabular}%
	\end{sidewaystable}
	
	\begin{sidewaystable}
		\caption{Size results for DGP.S5--DGP.S7 for the procedure based on $\Delta_a$.}\label{tab:S2}
		\centering\scriptsize
		\addtolength{\tabcolsep}{-2pt}
		\begin{tabular}{cccccccccccccccccccccccccccccccccccccc}
			\toprule
			&&& \multicolumn{5}{c}{DGP.S5} & \multicolumn{5}{c}{DGP.S6} & \multicolumn{5}{c}{DGP.S7} \\
			\cmidrule(lr){4-8} \cmidrule(lr){9-13} \cmidrule(lr){14-18}
			$m$ & $T$ & $L$
			& \multicolumn{1}{c}{$\Delta_{0.1}$} & \multicolumn{1}{c}{$\Delta_{0.5}$} & \multicolumn{1}{c}{$\Delta_{1}$} & \multicolumn{1}{c}{$\Delta_{1.5}$} & \multicolumn{1}{c}{$\Delta_{5}$}
			& \multicolumn{1}{c}{$\Delta_{0.1}$} & \multicolumn{1}{c}{$\Delta_{0.5}$} & \multicolumn{1}{c}{$\Delta_{1}$} & \multicolumn{1}{c}{$\Delta_{1.5}$} & \multicolumn{1}{c}{$\Delta_{5}$}
			& \multicolumn{1}{c}{$\Delta_{0.1}$} & \multicolumn{1}{c}{$\Delta_{0.5}$} & \multicolumn{1}{c}{$\Delta_{1}$} & \multicolumn{1}{c}{$\Delta_{1.5}$} & \multicolumn{1}{c}{$\Delta_{5}$} \\
			\midrule
			1 & 100 & 1 & 9.9 & 9.1 & 6.2 & 5.3 & 4.4 & 6.8 & 7.5 & 7.3 & 6.6 & 7.3 & 4.3 & 2.9 & 2.8 & 2.7 & 2.9 \\
			& & 2 & 9.1 & 9.1 & 8.7 & 7.7 & 4.9 & 7.4 & 8.3 & 9.2 & 8.8 & 8.0 & 4.4 & 3.7 & 3.4 & 3.4 & 3.1 \\
			& & 3 & 8.4 & 7.7 & 6.8 & 6.4 & 5.2 & 7.0 & 10 & 8.5 & 8.1 & 9.0 & 3.3 & 5.1 & 5.2 & 4.4 & 2.8 \\ 
			\cmidrule{2-18}
			& 200 & 1 & 6.3 & 7.7 & 7.2 & 6.4 & 5.3 & 8.0 & 8.8 & 8.6 & 8.8 & 7.8 & 5.0 & 4.9 & 4.2 & 3.9 & 2.9 \\
			& & 2 & 8.2 & 7.2 & 7.1 & 6.2 & 6.2 & 6.8 & 6.5 & 6.7 & 6.1 & 6.8 & 5.3 & 5.7 & 5.4 & 5.9 & 5.4 \\
			& & 3 & 11.3 & 9.4 & 10 & 8.0 & 6.4 & 6.8 & 6.5 & 7.4 & 7.5 & 7.5 & 4.3 & 4.2 & 3.5 & 3.1 & 2.1 \\ 
			\cmidrule{2-18}
			& 300 & 1 & 11.8 & 10.2 & 9.4 & 9.6 & 7.3 & 6.8 & 6.5 & 7.4 & 7.5 & 7.5 & 4.3 & 3.6 & 4.2 & 3.6 & 3.2 \\
			& & 2 & 7.5 & 7.6 & 7.2 & 6.1 & 4.5 & 6.7 & 8.4 & 8.1 & 8.2 & 7.8 & 3.2 & 4.2 & 4.1 & 4.0 & 3.5 \\
			& & 3 & 10.8 & 8.0 & 8.0 & 7.6 & 6.3 & 8.2 & 8.8 & 10.0 & 9.2 & 8.1 & 4.6 & 4.6 & 4.7 & 4.8 & 4.3 \\ 
			\midrule
			2 & 100 & 1 & 9.9 & 6.9 & 7.1 & 6.7 & 4.5 & 6.6 & 7.3 & 6.9 & 6.7 & 7.7 & 6.6 & 8.2 & 7.3 & 6.8 & 7.2 \\
			& & 2 & 8.3 & 8.4 & 7.2 & 7.3 & 5.5 & 6.8 & 8.7 & 8.5 & 8.8 & 7.1 & 5.6 & 4.0 & 3.2 & 4.0 & 3.8 \\
			& & 3 & 9.6 & 6.7 & 6.4 & 5.9 & 5.6 & 5.5 & 6.8 & 7.5 & 7.6 & 7.7 & 3.7 & 4.9 & 5.2 & 4.1 & 3.1 \\ 
			\cmidrule{2-18}
			& 200 & 1 & 7.9 & 6.9 & 8.0 & 7.0 & 5.3 & 6.6 & 8.2 & 7.3 & 6.8 & 7.2 & 5.6 & 4.9 & 4.9 & 4.5 & 3.1 \\
			& & 2 & 8.7 & 5.7 & 6.4 & 6.4 & 5.1 & 4.9 & 7.2 & 7.5 & 7.3 & 9.5 & 6.1 & 6.1 & 5.7 & 4.9 & 5.0 \\
			& & 3 & 10.8 & 8.4 & 8.2 & 8.9 & 6.0 & 6.6 & 7.0 & 6.7 & 6.7 & 7.1 & 4.3 & 3.0 & 3.0 & 2.2 & 2.8 \\ 
			\cmidrule{2-18}
			& 300 & 1 & 10.7 & 9.9 & 9.0 & 9.0 & 7.8 & 6.5 & 7.2 & 6.9 & 6.5 & 6.6 & 4.6 & 4.0 & 4.2 & 4.4 & 3.7 \\
			& & 2 & 8.5 & 7.1 & 6.7 & 5.3 & 4.4 & 7.0 & 7.8 & 7.5 & 7.2 & 7.0 & 3.2 & 3.8 & 3.8 & 4.2 & 3.3 \\
			& & 3 & 9.0 & 8.0 & 6.4 & 6.3 & 5.9 & 7.6 & 8.6 & 8.3 & 8.2 & 8.4 & 5.0 & 5.7 & 5.6 & 5.8 & 4.8 \\ 
			\midrule
			4 & 100 & 1 & 9.2 & 7.9 & 7.9 & 7.0 & 5.5 & 4.6 & 7.2 & 7.2 & 7.4 & 6.8 & 4.8 & 4.7 & 4.3 & 4.3 & 2.5 \\
			& & 2 & 9.4 & 7.8 & 6.0 & 6.7 & 6.1 & 5.7 & 8.2 & 7.5 & 7.5 & 6.8 & 4.9 & 4.7 & 4.2 & 3.7 & 2.9 \\
			& & 3 & 10.0 & 7.4 & 7.1 & 6.7 & 5.7 & 4.7 & 6.7 & 6.5 & 6.7 & 6.7 & 4.9 & 4.8 & 4.1 & 3.8 & 3.0 \\ 
			\cmidrule{2-18}
			& 200 & 1 & 8.4 & 7.4 & 6.5 & 6.9 & 5.0 & 5.1 & 6.6 & 6.8 & 7.1 & 6.8 & 4.8 & 5.0 & 4.9 & 4.0 & 2.9 \\
			& & 2 & 9.8 & 7.0 & 6.0 & 5.5 & 4.9 & 6.1 & 7.9 & 7.7 & 7.7 & 7.9 & 4.8 & 4.5 & 4.0 & 4.1 & 3.0 \\
			& & 3 & 12.6 & 9.3 & 7.7 & 7.5 & 6.7 & 4.6 & 6.1 & 6.6 & 6.9 & 7.2 & 4.4 & 3.7 & 3.1 & 3.1 & 2.7 \\ 
			\cmidrule{2-18}
			& 300 & 1 & 13.5 & 11.5 & 9.7 & 8.3 & 9.3 & 6.1 & 6.5 & 6.4 & 6.3 & 6.7 & 5.1 & 5.0 & 4.7 & 4.4 & 3.8 \\
			& & 2 & 12.0 & 8.5 & 7.8 & 7.3 & 5.2 & 5.9 & 8.0 & 8.4 & 8.2 & 7.5 & 4.7 & 4.8 & 4.7 & 4.5 & 3.2 \\
			& & 3 & 12.8 & 8.2 & 7.7 & 6.0 & 5.3 & 6.5 & 6.9 & 6.7 & 7.1 & 7.6 & 5.1 & 4.7 & 4.3 & 4.0 & 3.6 \\
			\bottomrule
		\end{tabular}%
	\end{sidewaystable}
	
	\begin{sidewaystable}
		\caption{Size results for DGP.S1--DGP.S7 for the procedures based on $\Delta^{(K)}$, $\Delta^{(L)}$ and $\Delta^{(R)}$.}\label{tab:Sother}
		\centering\scriptsize
		\addtolength{\tabcolsep}{-2pt}
		\begin{tabular}{cccccccccccccccccccccccccccccccccccccc}
			\toprule
			&& \multicolumn{3}{c}{DGP.S1} & \multicolumn{3}{c}{DGP.S2} & \multicolumn{3}{c}{DGP.S3} & \multicolumn{3}{c}{DGP.S4} & \multicolumn{3}{c}{DGP.S5} & \multicolumn{3}{c}{DGP.S6} & \multicolumn{3}{c}{DGP.S7} \\
			\cmidrule(lr){3-5} \cmidrule(lr){6-8} \cmidrule(lr){9-11} \cmidrule(lr){12-14} \cmidrule(lr){15-17} \cmidrule(lr){18-20} \cmidrule(lr){21-23}
			$T$ & $L$ &
			{$\Delta^{(K)}$} & {$\Delta^{(L)}$} & {$\Delta^{(R)}$} & 
			{$\Delta^{(K)}$} & {$\Delta^{(L)}$} & {$\Delta^{(R)}$} & 
			{$\Delta^{(K)}$} & {$\Delta^{(L)}$} & {$\Delta^{(R)}$} & 
			{$\Delta^{(K)}$} & {$\Delta^{(L)}$} & {$\Delta^{(R)}$} & 
			{$\Delta^{(K)}$} & {$\Delta^{(L)}$} & {$\Delta^{(R)}$} & 
			{$\Delta^{(K)}$} & {$\Delta^{(L)}$} & {$\Delta^{(R)}$} & 
			{$\Delta^{(K)}$} & {$\Delta^{(L)}$} & {$\Delta^{(R)}$} \\
			\midrule
			100 & 1 &  4.2 &  3.9 &  5.0 &  9.5 &  6.6 & 12.7 &  4.6 &  7.4 &  7.1 &  6.7 &  9.0 &  6.9 &  5.6 &  7.3 & 11.0 &  4.2 &  5.9 &  6.8 &  4.4 &  3.2 &  5.5\\
			& 2 &  5.5 &  4.4 &  4.1 &  8.3 &  4.9 & 15.4 &  4.0 &  6.2 &  7.4 &  5.3 &  8.6 &  6.7 &  6.6 &  7.4 &  8.6 &  5.9 &  2.0 &  3.9 &  5.4 &  4.9 &  5.0\\
			& 3 &  6.5 &  3.3 &  5.0 & 10.2 &  4.7 & 16.4 &  5.4 &  5.9 &  4.5 &  6.5 &  9.9 &  6.9 &  5.7 &  7.0 &  8.5 &  3.8 &  2.6 &  3.7 &  3.2 &  3.1 &  4.6\\ \midrule
			200 & 1 &  5.3 &  4.5 &  5.5 &  7.8 &  6.0 & 12.6 &  4.7 &  5.4 &  4.2 &  6.5 & 10.3 &  7.6 &  8.0 &  6.5 &  8.5 &  4.1 &  6.2 &  4.3 &  4.7 &  4.2 &  5.4\\
			& 2 &  5.0 &  3.8 &  5.6 &  6.3 &  4.2 & 10.4 &  6.7 &  5.7 &  5.8 &  4.1 &  5.6 &  8.0 &  8.2 &  7.3 &  8.3 &  5.8 &  2.4 &  6.8 &  4.1 &  3.1 &  5.1\\
			& 3 &  3.7 &  2.7 &  4.4 &  7.0 &  3.6 & 14.1 &  4.8 &  4.5 &  8.0 &  5.7 &  9.0 &  9.5 &  8.3 &  5.2 &  9.0 &  4.0 &  4.0 &  4.0 &  0.0 &  3.1 &  6.3\\ \midrule
			300 & 1 &  7.4 &  2.5 &  3.9 &  6.8 &  3.3 & 10.1 &  6.1 &  7.9 &  6.3 &  7.5 &  9.6 &  7.6 &  7.4 &  8.6 &  9.0 &  6.4 &  5.2 &  4.5 &  2.7 &  4.9 &  3.5\\
			& 2 &  5.4 &  3.4 &  5.1 &  5.9 &  4.1 & 10.6 &  6.2 &  4.7 &  6.7 &  6.0 &  6.8 & 10.6 &  7.9 &  8.9 &  7.8 &  4.0 &  3.8 &  7.4 &  0.0 &  3.3 &  3.1\\
			& 3 &  6.5 &  2.9 &  4.1 &  8.7 &  3.2 &  9.6 &  6.6 &  3.7 &  7.0 &  6.4 &  7.7 &  9.9 &  6.7 &  6.4 &  8.1 &  4.4 &  2.9 &  7.8 &  4.4 &  2.7 &  4.2\\
			\bottomrule
		\end{tabular}%
	\end{sidewaystable}
	
	\begin{sidewaystable}
		\caption{Power results for DGP.P1--DGP.P3 for the procedure based on $\Delta_a$.}\label{tab:P1}
		\centering\scriptsize
		\addtolength{\tabcolsep}{-2pt}
		\begin{tabular}{ccccccccccccccccccccc}
			\toprule
			&&& \multicolumn{5}{c}{DGP.P1} & \multicolumn{5}{c}{DGP.P2} & \multicolumn{5}{c}{DGP.P3} \\
			\cmidrule(lr){4-8} \cmidrule(lr){9-13} \cmidrule(lr){14-18}
			$m$ & $T$ & $L$
			& \multicolumn{1}{c}{$\Delta_{0.1}$} & \multicolumn{1}{c}{$\Delta_{0.5}$} & \multicolumn{1}{c}{$\Delta_{1}$} & \multicolumn{1}{c}{$\Delta_{1.5}$} & \multicolumn{1}{c}{$\Delta_{5}$}
			& \multicolumn{1}{c}{$\Delta_{0.1}$} & \multicolumn{1}{c}{$\Delta_{0.5}$} & \multicolumn{1}{c}{$\Delta_{1}$} & \multicolumn{1}{c}{$\Delta_{1.5}$} & \multicolumn{1}{c}{$\Delta_{5}$}
			& \multicolumn{1}{c}{$\Delta_{0.1}$} & \multicolumn{1}{c}{$\Delta_{0.5}$} & \multicolumn{1}{c}{$\Delta_{1}$} & \multicolumn{1}{c}{$\Delta_{1.5}$} & \multicolumn{1}{c}{$\Delta_{5}$} \\
			\midrule
			1 & 100 & 1 & 74.6 & 84.2 & 87.0 & 87.7 & 88.9 & 64.3 & 74.4 & 76.6 & 76.1 & 62.7 & 68.9 & 49.7 & 33.6 & 23.2 & 7.9 \\
			& & 2 & 78.8 & 87.8 & 89.2 & 90.1 & 90.5 & 70.3 & 83.6 & 86.4 & 86.1 & 74.1 & 76.7 & 58.4 & 38.2 & 28.1 & 10.3 \\
			& & 3 & 83.4 & 89.4 & 90.6 & 91.3 & 91.5 & 73.9 & 84.7 & 86.4 & 86.7 & 70.0 & 81.2 & 65.8 & 45.4 & 31.3 & 8.9 \\ 
			\cmidrule{2-18}
			& 200 & 1 & 87.2 & 93.9 & 95.8 & 96.0 & 96.7 & 78.9 & 89.9 & 90.7 & 91.0 & 85.6 & 85.5 & 73.7 & 55.4 & 42.2 & 14.0 \\
			& & 2 & 90.7 & 95.7 & 96.7 & 96.7 & 97.4 & 85.0 & 94.4 & 96.1 & 96.2 & 93.8 & 91.6 & 78.9 & 62.7 & 50.5 & 14.2 \\
			& & 3 & 93.3 & 97.2 & 98.0 & 98.5 & 98.6 & 85.6 & 94.9 & 95.6 & 95.8 & 92.9 & 91.8 & 80.1 & 67.2 & 54.8 & 12.0 \\ 
			\cmidrule{2-18}
			& 300 & 1 & 93.0 & 97.0 & 98.0 & 98.2 & 98.4 & 86.2 & 94.4 & 95.3 & 95.6 & 92.4 & 93.1 & 79.7 & 65.9 & 51.2 & 14.8 \\
			& & 2 & 95.8 & 98.2 & 99.0 & 99.2 & 99.5 & 92.2 & 96.9 & 98.3 & 98.5 & 96.3 & 95.1 & 84.8 & 75.2 & 62.2 & 15.1 \\
			& & 3 & 97.0 & 98.5 & 99.0 & 99.1 & 99.3 & 92.6 & 98.1 & 99.0 & 99.0 & 98.2 & 97.3 & 89.9 & 81.2 & 71.5 & 21.2 \\ 
			\midrule
			2 & 100 & 1 & 66.9 & 81.1 & 85.0 & 87.1 & 88.7 & 54.9 & 73.0 & 76.8 & 77.4 & 70.2 & 74.0 & 52.7 & 35.0 & 24.8 & 7.7 \\
			& & 2 & 73.6 & 86.6 & 88.7 & 89.2 & 90.3 & 64.3 & 83.3 & 87.3 & 88.1 & 82.2 & 80.4 & 61.6 & 43.2 & 32.8 & 10.5 \\
			& & 3 & 77.7 & 88.8 & 90.2 & 90.6 & 91.3 & 67.1 & 84.6 & 87.5 & 87.9 & 78.1 & 82.7 & 67.4 & 48.4 & 33.9 & 9.4 \\ 
			\cmidrule{2-18}
			& 200 & 1 & 81.4 & 92.0 & 94.7 & 95.4 & 96.4 & 74.3 & 87.3 & 90.8 & 91.4 & 88.2 & 88.1 & 73.2 & 55.8 & 41.6 & 13.0 \\
			& & 2 & 88.2 & 95.1 & 96.2 & 96.6 & 97.2 & 81.8 & 93.2 & 96.3 & 96.8 & 96.0 & 93.5 & 80.1 & 65.6 & 53.4 & 13.9 \\
			& & 3 & 91.0 & 96.8 & 97.7 & 98.0 & 98.7 & 81.5 & 94.5 & 96.2 & 96.5 & 95.2 & 93.8 & 82.1 & 70.0 & 56.8 & 12.7 \\ 
			\cmidrule{2-18}
			& 300 & 1 & 89.1 & 95.9 & 97.4 & 97.7 & 98.4 & 81.8 & 93.6 & 95.9 & 96.1 & 94.9 & 94.8 & 92.0 & 68.3 & 54.6 & 14.9 \\
			& & 2 & 92.6 & 97.6 & 98.6 & 99.2 & 99.6 & 87.9 & 96.3 & 98.3 & 98.6 & 97.9 & 96.7 & 87.5 & 76.1 & 63.5 & 16.4 \\
			& & 3 & 94.7 & 98.0 & 98.7 & 99.0 & 99.4 & 89.0 & 97.0 & 98.9 & 99.2 & 99.0 & 98.1 & 90.8 & 92.0 & 72.9 & 22.9 \\ 
			\midrule
			4 & 100 & 1 & 45.3 & 73.8 & 81.3 & 83.7 & 87.5 & 39.9 & 67.6 & 76.5 & 79.9 & 78.4 & 74.4 & 54.7 & 36.6 & 26.8 & 8.9 \\
			& & 2 & 52.8 & 79.4 & 85.3 & 87.9 & 90.5 & 50.5 & 76.7 & 83.7 & 86.3 & 86.4 & 84.0 & 64.6 & 46.6 & 33.9 & 11.3 \\
			& & 3 & 56.9 & 81.0 & 86.4 & 88.5 & 91.2 & 53.6 & 79.2 & 86.0 & 88.5 & 88.1 & 86.1 & 68.6 & 50.6 & 37.4 & 10.6 \\ 
			\cmidrule{2-18}
			& 200 & 1 & 65.0 & 86.6 & 91.5 & 93.1 & 94.9 & 58.9 & 81.4 & 88.7 & 91.1 & 92.3 & 89.5 & 73.8 & 57.8 & 44.4 & 13.2 \\
			& & 2 & 73.4 & 90.7 & 94.8 & 95.8 & 97.2 & 66.3 & 87.0 & 92.7 & 94.6 & 95.5 & 95.1 & 81.5 & 67.6 & 55.4 & 16.0 \\
			& & 3 & 76.0 & 92.7 & 95.8 & 97.0 & 98.2 & 71.0 & 90.3 & 95.1 & 96.7 & 97.5 & 95.8 & 83.3 & 71.2 & 59.1 & 17.1 \\ 
			\cmidrule{2-18}
			& 300 & 1 & 74.7 & 92.0 & 95.1 & 96.3 & 97.7 & 69.1 & 88.4 & 94.3 & 95.9 & 96.9 & 95.5 & 83.3 & 69.8 & 57.7 & 16.2 \\
			& & 2 & 81.6 & 95.5 & 97.7 & 98.2 & 99.2 & 74.9 & 93.6 & 97.2 & 98.4 & 99.1 & 98.1 & 89.4 & 78.5 & 67.4 & 19.2 \\
			& & 3 & 83.7 & 96.9 & 98.5 & 99.2 & 3.0 & 78.2 & 94.2 & 97.9 & 98.9 & 99.5 & 98.7 & 91.1 & 82.0 & 73.0 & 26.3 \\
			\bottomrule
		\end{tabular}
	\end{sidewaystable}
	
	\begin{sidewaystable}
		\caption{Power results for DGP.P4--DGP.P5 for the procedure based on $\Delta_a$.}\label{tab:P2}
		\centering\scriptsize
		\addtolength{\tabcolsep}{-2pt}
		\begin{tabular}{cccccccccccccccccccccccccccccccccccccc}
			\toprule
			&&& \multicolumn{5}{c}{DGP.P4} & \multicolumn{5}{c}{DGP.P5} \\
			\cmidrule(lr){4-8} \cmidrule(lr){9-13}
			$m$ & $T$ & $L$
			& \multicolumn{1}{c}{$\Delta_{0.1}$} & \multicolumn{1}{c}{$\Delta_{0.5}$} & \multicolumn{1}{c}{$\Delta_{1}$} & \multicolumn{1}{c}{$\Delta_{1.5}$} & \multicolumn{1}{c}{$\Delta_{5}$}
			& \multicolumn{1}{c}{$\Delta_{0.1}$} & \multicolumn{1}{c}{$\Delta_{0.5}$} & \multicolumn{1}{c}{$\Delta_{1}$} & \multicolumn{1}{c}{$\Delta_{1.5}$} & \multicolumn{1}{c}{$\Delta_{5}$} \\
			\midrule
			1 & 100 & 1 & 17.5 & 21.8 & 20.3 & 15.9 & 11.7 & 32.1 & 42.5 & 45.7 & 45.9 & 42.2 \\
			& & 2 & 32.0 & 51.6 & 49.8 & 44.5 & 19.2 & 38.7 & 48.1 & 47.0 & 46.9 & 38.4 \\
			& & 3 & 47.6 & 59.6 & 55.7 & 45.7 & 16.4 & 46.5 & 47.4 & 45.8 & 46.3 & 37.5 \\ 
			\cmidrule{2-13}
			& 200 & 1 & 36.5 & 56.8 & 56.4 & 49.5 & 19.6 & 56.0 & 53.8 & 61.5 & 61.9 & 58.8 \\
			& & 2 & 66.7 & 84.8 & 85.5 & 79.5 & 32.8 & 6.8 & 71.3 & 69.2 & 70.0 & 64.2 \\
			& & 3 & 76.7 & 86.0 & 85.7 & 81.2 & 30.9 & 71.8 & 71.5 & 70.1 & 68.5 & 58.5 \\ 
			\cmidrule{2-13}
			& 300 & 1 & 52.0 & 77.7 & 77.6 & 71.7 & 32.8 & 66.5 & 72.7 & 73.9 & 74.7 & 72.0 \\
			& & 2 & 89.8 & 96.4 & 96.9 & 95.5 & 64.5 & 74.2 & 78.4 & 76.1 & 74.9 & 70.9 \\
			& & 3 & 96.0 & 99.0 & 99.2 & 98.9 & 75.5 & 81.8 & 83.7 & 83.4 & 82.5 & 75.7 \\ 
			\midrule
			2 & 100 & 1 & 11.8 & 24.3 & 28.0 & 23.1 & 13.6 & 37.5 & 52.8 & 57.4 & 57.8 & 54.0 \\
			& & 2 & 30.3 & 57.2 & 58.8 & 53.9 & 26.3 & 46.8 & 55.9 & 57.8 & 57.0 & 51.6 \\
			& & 3 & 41.7 & 60.7 & 63.8 & 57.5 & 20.9 & 51.9 & 57.1 & 55.5 & 55.5 & 51.3 \\ 
			\cmidrule{2-13}
			& 200 & 1 & 32.5 & 58.8 & 61.5 & 57.0 & 28.5 & 59.8 & 69.9 & 72.5 & 72.5 & 70.4 \\
			& & 2 & 62.6 & 86.2 & 89.3 & 87.9 & 52.7 & 69.7 & 77.6 & 77.9 & 77.4 & 75.6 \\
			& & 3 & 70.5 & 90.3 & 91.7 & 91.0 & 50.7 & 73.8 & 78.2 & 77.6 & 76.7 & 71.1 \\ 
			\cmidrule{2-13}
			& 300 & 1 & 48.0 & 80.5 & 84.3 & 82.3 & 51.7 & 67.7 & 78.6 & 80.0 & 80.1 & 80.3 \\
			& & 2 & 83.0 & 96.9 & 97.9 & 97.8 & 84.7 & 77.4 & 83.7 & 83.6 & 81.7 & 79.0 \\
			& & 3 & 92.8 & 99.2 & 99.4 & 99.4 & 91.7 & 81.7 & 88.9 & 88.1 & 88.1 & 84.3 \\ 
			\midrule
			4 & 100 & 1 & 6.6 & 25.4 & 33.0 & 34.9 & 22.3 & 33.2 & 54.6 & 59.7 & 61.1 & 61.9 \\
			& & 2 & 17.7 & 49.8 & 58.9 & 61.2 & 36.7 & 46.4 & 63.4 & 66.0 & 67.9 & 67.4 \\
			& & 3 & 24.8 & 55.7 & 65.1 & 66.8 & 38.1 & 51.2 & 65.5 & 67.4 & 68.0 & 67.7 \\ 
			\cmidrule{2-13}
			& 200 & 1 & 16.3 & 49.1 & 62.2 & 65.3 & 42.6 & 50.8 & 72.3 & 77.0 & 78.9 & 79.2 \\
			& & 2 & 42.8 & 82.9 & 91.1 & 91.6 & 73.8 & 64.3 & 79.0 & 82.0 & 82.2 & 81.1 \\
			& & 3 & 56.1 & 90.2 & 95.1 & 95.7 & 83.7 & 68.7 & 81.5 & 83.1 & 83.2 & 81.1 \\ 
			\cmidrule{2-13}
			& 300 & 1 & 27.1 & 73.4 & 86.4 & 88.8 & 72.7 & 62.3 & 80.5 & 83.5 & 84.2 & 85.0 \\
			& & 2 & 69.0 & 95.8 & 98.8 & 99.0 & 95.0 & 72.0 & 85.9 & 87.7 & 87.8 & 85.9 \\
			& & 3 & 78.3 & 98.4 & 99.4 & 99.5 & 97.6 & 77.6 & 88.5 & 89.4 & 89.2 & 86.5 \\
			\bottomrule
		\end{tabular}%
	\end{sidewaystable}
	
	\begin{sidewaystable}
		\caption{Power results for DGP.P1--DGP.P5 for the procedures based on $\Delta^{(K)}$, $\Delta^{(L)}$ and $\Delta^{(R)}$}\label{tab:Pother}
		\centering\scriptsize
		\addtolength{\tabcolsep}{-2pt}
		\begin{tabular}{cccccccccccccccccccccccccccccccccccccc}
			\toprule
			&& \multicolumn{3}{c}{DGP.P1} & \multicolumn{3}{c}{DGP.P2} & \multicolumn{3}{c}{DGP.P3} & \multicolumn{3}{c}{DGP.P4} & \multicolumn{3}{c}{DGP.P5} \\
			\cmidrule(lr){3-5} \cmidrule(lr){6-8} \cmidrule(lr){9-11} \cmidrule(lr){12-14} \cmidrule(lr){15-17} 
			$T$ & $L$ &
			{$\Delta^{(K)}$} & {$\Delta^{(L)}$} & {$\Delta^{(R)}$} & 
			{$\Delta^{(K)}$} & {$\Delta^{(L)}$} & {$\Delta^{(R)}$} & 
			{$\Delta^{(K)}$} & {$\Delta^{(L)}$} & {$\Delta^{(R)}$} & 
			{$\Delta^{(K)}$} & {$\Delta^{(L)}$} & {$\Delta^{(R)}$} & 
			{$\Delta^{(K)}$} & {$\Delta^{(L)}$} & {$\Delta^{(R)}$} \\
			\midrule
			100 & 1 &  72.2 &  69.6 &  98.0 &  64.3 &  66.7 &  29.8 &  40.5 &  27.0 &  22.3 &  34.5 &  30.0 &  15.3 &  27.5 &  26.4 &  15.2\\
			& 2 &  76.5 &  84.6 & 100.0 &  68.0 &  77.8 &  39.6 &  49.1 &  38.0 &  41.8 &  39.6 &  31.5 &  14.8 &  35.1 &  35.4 &  17.9\\
			& 3 &  83.3 &  88.5 & 100.0 &  72.4 &  80.4 &  54.8 &  54.4 &  42.7 &  59.0 &  47.2 &  20.2 &  14.1 &  41.1 &  37.1 &  20.5\\ \midrule
			200 & 1 &  83.2 &  82.4 &  99.9 &  71.0 &  82.2 &  57.0 &  64.7 &  47.8 &  61.0 &  46.3 &  48.9 &  19.5 &  44.8 &  47.2 &  18.8\\
			& 2 &  87.8 &  93.4 & 100.0 &  80.9 &  92.6 &  82.0 &  70.6 &  56.9 &  90.1 &  67.5 &  38.7 &  15.9 &  54.1 &  43.6 &  25.5\\
			& 3 &  88.3 &  96.8 & 100.0 &  81.8 &  89.5 &  92.2 &  73.7 &  55.2 &  97.2 &  58.1 &  24.3 &  18.6 &  49.0 &  43.8 &  19.5\\ \midrule
			300 & 1 &  87.4 &  87.4 & 100.0 &  81.6 &  89.8 &  71.7 &  74.8 &  56.2 &  82.5 &  68.4 &  57.5 &  21.4 &  54.0 &  53.6 &  25.3\\
			& 2 &  91.0 &  96.5 & 100.0 &  85.1 &  95.9 &  96.4 &  79.8 &  67.1 &  97.6 &  81.9 &  43.5 &  22.8 &  61.1 &  60.8 &  30.5\\
			& 3 &  91.4 &  96.8 & 100.0 &  88.6 &  91.4 &  99.5 &  83.9 &  70.1 &  99.8 &  79.7 &  27.4 &  21.3 &  64.7 &  51.7 &  36.7\\
			\bottomrule
		\end{tabular}%
	\end{sidewaystable}
	
	\section{Real-data application}
	In this section we provide an application to financial market data
	to demonstrate the potential of the monitoring procedure to detect breaks in stationarity.
	We identified all constituent stocks of the S\&P500 stock index which
	exhibited no structural change in weekly returns over the six-year period spanning from
	1~January~2002 to 31~December~2007 (\emph{training period}).
	We then ran the monitoring procedure on weekly returns starting on 1~January~2008 and continued
	the monitoring for eight and a half years (until 30~June~2016) or until the first break in stationarity occurred,
	whichever came first.
	All returns were based on weekly closing prices adjusted for dividends and splits.
	The data that support the findings of this study are openly available for download via the Alpha Vantage API \citep{alphavantage}.
	
	The initial dataset consisted of 311 companies which remained components of the S\&P500 index
	throughout the total period under consideration.
	To identify which stocks exhibited no significant break during the training period,
	we employed the test for strict stationarity recently proposed by \citet{hong2016}.
	Based on this test 71 of the stocks were identified as having no structural breaks
	in stationarity of dimension $m=1$ (i.e.\ no change in marginal distribution) over the training period
	(at the $5\%$ level of significance).
	The 71 identified stocks are listed in Table~\ref{tab:stocks1}.
	
	\begin{table}
		\caption{Constituent stocks of the S\&P500 index which, by the test of \citet{hong2016}, exhibited no significant
			structural breaks (at $\alpha=5\%$) from 1~January~2002 to 31~December~2007.}\label{tab:stocks1}
		\centering\tiny
		\addtolength{\tabcolsep}{-1pt}
		\def\arraystretch{.7}
		\begin{tabular}{lcccclcccclcccc}
			\toprule
			Stock & $m$ & $\hat{p}$ & $\hat\tau$ & Date & Stock & $m$ & $\hat{p}$ & $\hat\tau$ & Date & Stock & $m$ & $\hat{p}$ & $\hat\tau$ & Date \\
			\cmidrule(lr){1-5} \cmidrule(lr){6-10} \cmidrule(lr){11-15}
			AA    & 1     & 0.000 & 52    &  08/12/29  &  DRI  & 1     & 0.000 & 45    &  08/11/10  &  MUR  & 1     & 0.001 & 62    &  09/03/09 \\
			& 2     & 0.000 & 54    &  09/01/12  &       & 2     & 0.001 & 45    &  08/11/10  &       & 2     & 0.000 & 62    &  09/03/09 \\
			& 4     & 0.000 & 58    &  09/02/09  &       & 4     & 0.004 & 52    &  08/12/29  &       & 4     & 0.001 & 62    &  09/03/09 \\
			AAPL  & 1     & 0.050 &  $\infty$  &  -    &  DVN  & 1     & 0.001 & 60    &  09/02/23  &  MYL  & 1     & 0.019 & 70    &  09/05/04 \\
			& 2     & 0.019 & 397   &  15/08/10  &       & 2     & 0.000 & 60    &  09/02/23  &       & 2     & 0.013 & 71    &  09/05/11 \\
			& 4     & 0.005 & 354   &  14/10/13  &       & 4     & 0.003 & 70    &  09/05/04  &       & 4     & 0.019 & 71    &  09/05/11 \\
			ADBE  & 1     & 0.256 &  $\infty$  &  -    &  EA   & 1     & 0.005 & 62    &  09/03/09  &  NUE  & 1     & 0.024 & 394   &  15/07/20 \\
			& 2     & 0.282 &  $\infty$  &  -    &       & 2     & 0.002 & 62    &  09/03/09  &       & 2     & 0.004 & 359   &  14/11/17 \\
			& 4     & 0.067 &  $\infty$  &  -    &       & 4     & 0.003 & 109   &  10/02/01  &       & 4     & 0.000 & 323   &  14/03/10 \\
			AET   & 1     & 0.000 & 54    &  09/01/12  &  EMN  & 1     & 0.000 & 60    &  09/02/23  &  PDCO  & 1     & 0.053 &  $\infty$  & - \\
			& 2     & 0.001 & 57    &  09/02/02  &       & 2     & 0.000 & 62    &  09/03/09  &       & 2     & 0.079 &  $\infty$  & - \\
			& 4     & 0.000 & 61    &  09/03/02  &       & 4     & 0.000 & 66    &  09/04/06  &       & 4     & 0.099 &  $\infty$  & - \\
			ALL   & 1     & 0.000 & 56    &  09/01/26  &  ESRX  & 1     & 0.098 &  $\infty$  &  -    &  PFE  & 1     & 0.059 &  $\infty$  & - \\
			& 2     & 0.000 & 57    &  09/02/02  &       & 2     & 0.048 & 442   &  16/06/20  &       & 2     & 0.048 & 396   &  15/08/03 \\
			& 4     & 0.000 & 60    &  09/02/23  &       & 4     & 0.011 & 391   &  15/06/29  &       & 4     & 0.028 & 379   &  15/04/06 \\
			AON   & 1     & 0.271 &  $\infty$  &  -    &  ETN  & 1     & 0.000 & 56    &  09/01/26  &  PG   & 1     & 0.004 & 61    &  09/03/02 \\
			& 2     & 0.159 &  $\infty$  &  -    &       & 2     & 0.000 & 56    &  09/01/26  &       & 2     & 0.001 & 61    &  09/03/02 \\
			& 4     & 0.057 &  $\infty$  &  -    &       & 4     & 0.000 & 57    &  09/02/02  &       & 4     & 0.007 & 69    &  09/04/27 \\
			APA   & 1     & 0.000 & 49    &  08/12/08  &  F    & 1     & 0.000 & 42    &  08/10/20  &  PLD  & 1     & 0.000 & 46    &  08/11/17 \\
			& 2     & 0.000 & 54    &  09/01/12  &       & 2     & 0.000 & 40    &  08/10/06  &       & 2     & 0.000 & 47    &  08/11/24 \\
			& 4     & 0.000 & 60    &  09/02/23  &       & 4     & 0.000 & 43    &  08/10/27  &       & 4     & 0.000 & 48    &  08/12/01 \\
			APC   & 1     & 0.000 & 48    &  08/12/01  &  FCX  & 1     & 0.000 & 52    &  08/12/29  &  PSA  & 1     & 0.000 & 50    &  08/12/15 \\
			& 2     & 0.000 & 49    &  08/12/08  &       & 2     & 0.001 & 57    &  09/02/02  &       & 2     & 0.000 & 53    &  09/01/05 \\
			& 4     & 0.000 & 51    &  08/12/22  &       & 4     & 0.001 & 74    &  09/06/01  &       & 4     & 0.000 & 55    &  09/01/20 \\
			BCR   & 1     & 0.565 &  $\infty$  &  -    &  FDX  & 1     & 0.000 & 53    &  09/01/05  &  R    & 1     & 0.000 & 50    &  08/12/15 \\
			& 2     & 0.561 &  $\infty$  &  -    &       & 2     & 0.000 & 49    &  08/12/08  &       & 2     & 0.000 & 50    &  08/12/15 \\
			& 4     & 0.679 &  $\infty$  &  -    &       & 4     & 0.000 & 51    &  08/12/22  &       & 4     & 0.000 & 52    &  08/12/29 \\
			BEN   & 1     & 0.000 & 47    &  08/11/24  &  GPC  & 1     & 0.001 & 59    &  09/02/17  &  ROK  & 1     & 0.000 & 54    &  09/01/12 \\
			& 2     & 0.000 & 47    &  08/11/24  &       & 2     & 0.000 & 53    &  09/01/05  &       & 2     & 0.000 & 54    &  09/01/12 \\
			& 4     & 0.000 & 54    &  09/01/12  &       & 4     & 0.000 & 56    &  09/01/26  &       & 4     & 0.000 & 57    &  09/02/02 \\
			BHI   & 1     & 0.000 & 73    &  09/05/26  &  GWW  & 1     & 0.065 &  $\infty$  &  -    &  STJ  & 1     & 0.429 &  $\infty$  & - \\
			& 2     & 0.000 & 73    &  09/05/26  &       & 2     & 0.064 &  $\infty$  &  -    &       & 2     & 0.414 &  $\infty$  & - \\
			& 4     & 0.000 & 79    &  09/07/06  &       & 4     & 0.107 &  $\infty$  &  -    &       & 4     & 0.462 &  $\infty$  & - \\
			CAG   & 1     & 0.118 &  $\infty$  &  -    &  HES  & 1     & 0.000 & 49    &  08/12/08  &  SYK  & 1     & 0.025 & 137   &  10/08/16 \\
			& 2     & 0.026 & 74    &  09/06/01  &       & 2     & 0.000 & 52    &  08/12/29  &       & 2     & 0.040 & 207   &  11/12/19 \\
			& 4     & 0.009 & 65    &  09/03/30  &       & 4     & 0.000 & 54    &  09/01/12  &       & 4     & 0.071 &  $\infty$  & - \\
			CAT   & 1     & 0.002 & 59    &  09/02/17  &  HOG  & 1     & 0.000 & 46    &  08/11/17  &  SYY  & 1     & 0.024 & 61    &  09/03/02 \\
			& 2     & 0.000 & 64    &  09/03/23  &       & 2     & 0.000 & 48    &  08/12/01  &       & 2     & 0.012 & 59    &  09/02/17 \\
			& 4     & 0.002 & 76    &  09/06/15  &       & 4     & 0.000 & 52    &  08/12/29  &       & 4     & 0.013 & 63    &  09/03/16 \\
			CHRW  & 1     & 0.130 &  $\infty$  &  -    &  HUM  & 1     & 0.011 & 57    &  09/02/02  &  UNH  & 1     & 0.000 & 37    &  08/09/15 \\
			& 2     & 0.099 &  $\infty$  &  -    &       & 2     & 0.010 & 57    &  09/02/02  &       & 2     & 0.000 & 37    &  08/09/15 \\
			& 4     & 0.052 &  $\infty$  &  -    &       & 4     & 0.007 & 61    &  09/03/02  &       & 4     & 0.000 & 39    &  08/09/29 \\
			CI    & 1     & 0.000 & 53    &  09/01/05  &  IP   & 1     & 0.000 & 47    &  08/11/24  &  UNM  & 1     & 0.019 & 70    &  09/05/04 \\
			& 2     & 0.000 & 54    &  09/01/12  &       & 2     & 0.000 & 47    &  08/11/24  &       & 2     & 0.020 & 71    &  09/05/11 \\
			& 4     & 0.000 & 60    &  09/02/23  &       & 4     & 0.000 & 50    &  08/12/15  &       & 4     & 0.039 & 84    &  09/08/10 \\
			CINF  & 1     & 0.000 & 41    &  08/10/13  &  ITW  & 1     & 0.000 & 74    &  09/06/01  &  UPS  & 1     & 0.000 & 49    &  08/12/08 \\
			& 2     & 0.000 & 40    &  08/10/06  &       & 2     & 0.000 & 61    &  09/03/02  &       & 2     & 0.000 & 44    &  08/11/03 \\
			& 4     & 0.000 & 43    &  08/10/27  &       & 4     & 0.000 & 71    &  09/05/11  &       & 4     & 0.000 & 44    &  08/11/03 \\
			CPB   & 1     & 0.036 & 68    &  09/04/20  &  KMB  & 1     & 0.016 & 62    &  09/03/09  &  UTX  & 1     & 0.000 & 59    &  09/02/17 \\
			& 2     & 0.019 & 61    &  09/03/02  &       & 2     & 0.017 & 68    &  09/04/20  &       & 2     & 0.000 & 61    &  09/03/02 \\
			& 4     & 0.010 & 62    &  09/03/09  &       & 4     & 0.020 & 71    &  09/05/11  &       & 4     & 0.003 & 66    &  09/04/06 \\
			CSX   & 1     & 0.000 & 57    &  09/02/02  &  LEG  & 1     & 0.000 & 53    &  09/01/05  &  VMC  & 1     & 0.000 & 47    &  08/11/24 \\
			& 2     & 0.000 & 58    &  09/02/09  &       & 2     & 0.000 & 53    &  09/01/05  &       & 2     & 0.000 & 48    &  08/12/01 \\
			& 4     & 0.001 & 60    &  09/02/23  &       & 4     & 0.000 & 55    &  09/01/20  &       & 4     & 0.000 & 46    &  08/11/17 \\
			DD    & 1     & 0.000 & 59    &  09/02/17  &  LM   & 1     & 0.000 & 45    &  08/11/10  &  WFM  & 1     & 0.002 & 47    &  08/11/24 \\
			& 2     & 0.000 & 61    &  09/03/02  &       & 2     & 0.000 & 49    &  08/12/08  &       & 2     & 0.000 & 47    &  08/11/24 \\
			& 4     & 0.000 & 76    &  09/06/15  &       & 4     & 0.000 & 52    &  08/12/29  &       & 4     & 0.001 & 56    &  09/01/26 \\
			DE    & 1     & 0.002 & 56    &  09/01/26  &  MMC  & 1     & 0.077 &  $\infty$  &  -    &  WHR  & 1     & 0.000 & 55    &  09/01/20 \\
			& 2     & 0.000 & 57    &  09/02/02  &       & 2     & 0.044 & 439   &  16/05/31  &       & 2     & 0.000 & 55    &  09/01/20 \\
			& 4     & 0.005 & 64    &  09/03/23  &       & 4     & 0.012 & 395   &  15/07/27  &       & 4     & 0.000 & 61    &  09/03/02 \\
			DHR   & 1     & 0.002 & 59    &  09/02/17  &  MMM  & 1     & 0.007 & 61    &  09/03/02  &  XL   & 1     & 0.000 & 27    &  08/07/07 \\
			& 2     & 0.000 & 60    &  09/02/23  &       & 2     & 0.004 & 63    &  09/03/16  &       & 2     & 0.000 & 27    &  08/07/07 \\
			& 4     & 0.001 & 69    &  09/04/27  &       & 4     & 0.004 & 74    &  09/06/01  &       & 4     & 0.000 & 28    &  08/07/14 \\
			DOV   & 1     & 0.000 & 54    &  09/01/12  &  MRK  & 1     & 0.026 & 62    &  09/03/09  &  XOM  & 1     & 0.068 &  $\infty$  & - \\
			& 2     & 0.000 & 58    &  09/02/09  &       & 2     & 0.035 & 63    &  09/03/16  &       & 2     & 0.060 &  $\infty$  & - \\
			& 4     & 0.000 & 63    &  09/03/16  &       & 4     & 0.047 & 101   &  09/12/07  &       & 4     & 0.049 & 84    &  09/08/10 \\
			DOW   & 1     & 0.000 & 59    &  09/02/17  &  MRO  & 1     & 0.001 & 52    &  08/12/29  &  YUM  & 1     & 0.019 & 61    &  09/03/02 \\
			& 2     & 0.000 & 60    &  09/02/23  &       & 2     & 0.000 & 57    &  09/02/02  &       & 2     & 0.008 & 63    &  09/03/16 \\
			& 4     & 0.000 & 63    &  09/03/16  &       & 4     & 0.000 & 55    &  09/01/20  &       & 4     & 0.047 & 97    &  09/11/09 \\
			\bottomrule
		\end{tabular}
	\end{table}
	
	The results of our monitoring procedure are also given in Table~\ref{tab:stocks1}.
	Based on the results from the Monte Carlo study, we chose the tuning parameter $a=1$ for practical application.
	For each stock the critical value $c_\alpha$ of the test was approximated
	by $B=1\,000$ bootstrap replications
	using the resampling procedure described in Section~\ref{sec:resampling.procedure}.
	The observed run length (in weeks, denoted by $\hat\tau$) of the monitoring procedure until a break in stationarity
	was identified is reported along with the corresponding calendar date.
	$\hat\tau=\infty$ indicates that no break in strict stationarity was encountered during monitoring.
	Additionally we also supply an approximate $p$-value, calculated as
	\[
	\hat{p} = \frac{1}{B} \sum_{b=1}^B \I \left( M^*_{(b)} \ge \max_{1 \le t \le LT} \Delta_{T,t} \right),
	\]
	where $M^*_{(b)}$ is as defined in Section~\ref{sec:resampling.procedure}.
	
	Most of the considered stock returns exhibit some degree of irregularity over the period running from 2008 to 2010.
	The obvious justification for these irregularities is the effects the global financial crisis
	had on markets during this period. Some sources claim that the crisis peaked in the fourth quarter of 2008 (see, e.g., \citealp{Ivashina15}).
	In agreement with this, our monitoring procedure identified significant structural breaks (at $\alpha=5\%$)
	in the returns of a large proportion of the stocks late in 2008 or early in 2009, as indicated in Table~\ref{tab:stocks1}.
	
	To further investigate the source of the breaks in stationarity we employ the methodology followed
	in \citet{white08} by means of a multi-quantile conditional autoregressive value at risk (MQ-CAViaR) model,
	which was originally proposed in its univariate form by \citet{engle04}. The objective is to model quantiles
	of the distribution of the return series, denoted by $r_t$, conditional on past observations
	$\{r_1,\ldots,r_{t-1}\}$, which may in turn be used to calculate measures of (conditional)
	distributional properties of $r_t$ over time. Specifically we are interested in the following
	measures, which are respectively measures of conditional spread (interquartile range, or IQR),
	conditional skewness (due to \citealp{bowley20}) and conditional kurtosis (due to \citealp{crow67}):
	\begin{equation*}\begin{split}
	K_{2,t} &= q_{4,t} - q_{2,t}, \\
	K_{3,t} &= \frac{q_{4,t} + q_{2,t} - 2q_{3,t}}{q_{4,t} - q_{2,t}}, \\
	K_{4,t} &= \frac{q_{5,t} - q_{1,t}}{q_{4,t} - q_{2,t}} - 2.91,
	\end{split}\end{equation*}
	where, for $(\theta_1,\ldots,\theta_5) = (0.025,0.25,0.5,0.75,0.975)$, each quantile $q_{j,t}$ is defined
	such that $\Prob(r_t \le q_{j,t} | \mathcal{F}_{t-1}) = \theta_j$. Here,
	$\mathcal{F}_{t-1}$ denotes the $\sigma$-algebra generated by $\{r_1,\ldots,r_{t-1}\}$.
	
\begin{figure}
		\centering
		\includegraphics[width=1.00\textwidth]{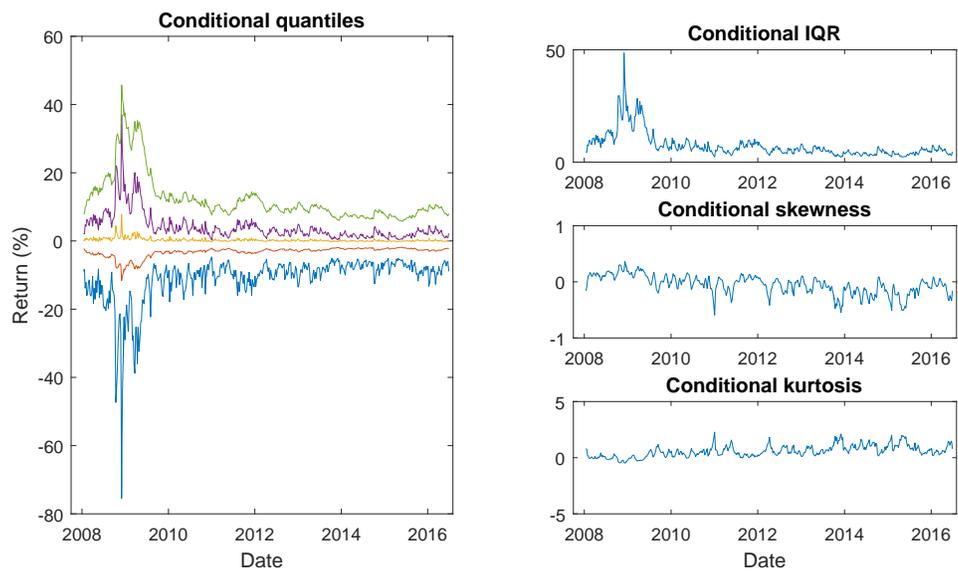}
		\caption{\emph{Left:}~0.025, 0.25, 0.5, 0.75 and 0.975 level conditional quantiles
			of weekly the returns on a share of the Ford Motor Company.
			\emph{Right:}~Conditional
			interquartile range and measures of skewness and kurtosis calculated from the conditional quantiles.}
		\label{fig:f_caviar}
	\end{figure}

	Because of the dramatic effect of the global economic crisis, many of the stock returns exhibit
	similar erratic behavior towards the end of 2008. As an example, we look at the returns of the Ford Motor Company.
	Graphical representations of this stock's conditional IQR, skewness and kurtosis,
	along with a graphical representation of the modeled quantiles, are given in Figure~\ref{fig:f_caviar}.
	The break in stationarity of returns seems to stem from
	a change in variability, which spiked just before the end of 2008.
	This was identified fairly accurately by our monitoring procedure, which terminated in the 10th month of 2008
	for $m=1,2,4$.
	Interestingly, the conditional skewness and kurtosis of these returns remained very stable
	during this period.
	
\begin{figure}
		\centering
		\includegraphics[width=1.00\textwidth]{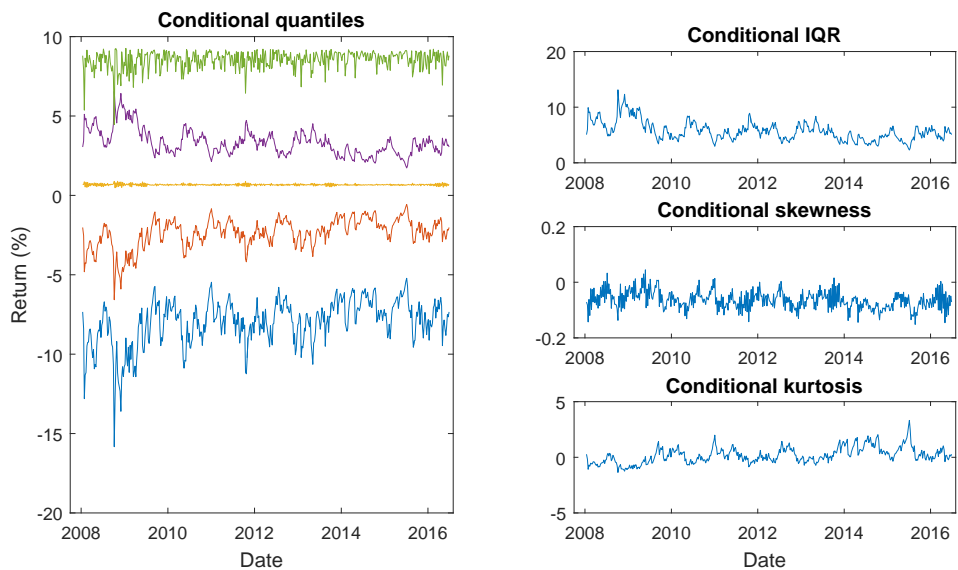}
		\caption{\emph{Left:}~0.025, 0.25, 0.5, 0.75 and 0.975 level conditional quantiles
			of weekly the returns on a share of Apple Inc.\
			\emph{Right:}~Conditional
			interquartile range and measures of skewness and kurtosis calculated from the conditional quantiles.}
		\label{fig:aapl_caviar}
	\end{figure}

\begin{figure}
		\centering
		\includegraphics[width=1.00\textwidth]{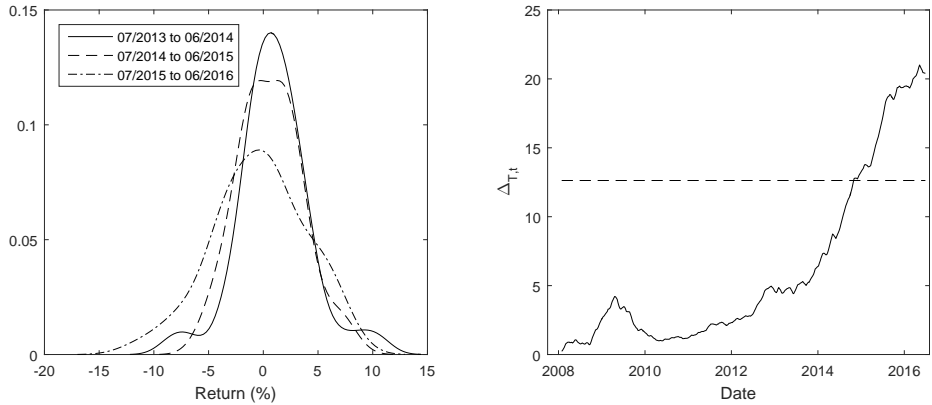}
		\caption{\emph{Left:}~Kernel density estimates for weekly returns of Apple Inc.\ stocks for three consecutive periods.
			\emph{Right:}~Value of the detector statistic over the monitoring period for $m=4$.
			The dashed line represents the estimated critical value~$\hat{c}_\alpha$.}
		\label{fig:aapl_dens}
	\end{figure}

	For some of the stocks the monitoring procedure rejected the null hypothesis at a much
	later stage. Consider, for example, the stock returns of Apple Inc.\ which were found to exhibit a break
	in stationarity in 2014 for $m=2$ and in 2015 for $m=4$. Looking at the MQ-CAViaR results
	presented in Figure~\ref{fig:aapl_caviar} no clear abrupt break in stationarity is apparent after 2009.
	However, close inspection of the conditional quantiles and kurtosis plots suggests a gradual change in
	distribution over time. This is confirmed by kernel density estimates of the weekly returns for three
	consecutive years, which are given in Figure~\ref{fig:aapl_dens}. The plot in the right hand side of
	Figure~\ref{fig:aapl_dens} confirms that our test statistic was able to detect these gradual structural changes.
	
	\section{Conclusion}
	We suggest a procedure for on-line monitoring of strict stationarity. The criteria involved are entirely model-free
	and therefore apply to arbitrary time series. In particular we employ a non-parametric estimate of the joint characteristic
	function of the underlying process and suggest to monitor an integrated functional of this estimate in order to capture structural
	breaks in the joint distribution of the observations. Asymptotic results include the limit null distribution of the detector statistic,
	as well as consistency under stationarity breaks of general nature. Since the limit null distribution depends on the stochastic structure
	of the series being sampled, a modification of the block bootstrap is used in order to actually implement the procedure. This bootstrap is
	in turn applied to simulated data, and shows satisfactory performance in terms of level and power. The same procedure is applied to real
	data from the financial market and the stationarity results are further scrutinized by incorporating the new MQ-CAViaR model.
	We close by reiterating that our results, although presented for scalar observations, they may readily be extended to multivariate time series.
	
	\subsection*{Acknowledgement}
	\sloppy
	The work of the third author was partially supported by OP RDE project No.\ CZ.02.2.69/0.0/0.0/16\_027/0008495, International Mobility of Researchers at Charles University.
	
	\section{Appendix}
	In this Appendix, we provide the proofs of Lemma 1 and Theorems 1
	and 2 in Section 4.
	\medskip\\ \indent {\bf Proof of Lemma 1}. Owing to Theorem 4 of \citet{kuelbs80}, we have that for some $0<\lambda<1/2$, with
	probability 1,
	\begin{eqnarray}\sum_{j=1}^k f(Y_j , u) -\sigma(u) {\cal B}(k) <<
	k^{1/2-\lambda}\ \ {\rm a.s.},\end{eqnarray} which results in the
	weak convergence of $Z_{T}(\theta, u)$ to $Z(\theta, u)$. In fact,
	the second statement of the theorem is established by applying
	Theorem 4 of \citet{kuelbs80} for stationary random vector
	process and Cram\'er-Wold's device. \hfill$\Box$
	\medskip\\ \indent
	{\bf Proof of Theorem 1}. We only consider the case that $\gamma=0$ because the other cases can be handled
	similarly.  Using  (1.12b) of \citet{rio13} with $p=2,
	q=r=4$ (see also Lemma 2.1 of \citet{kuelbs80}, it is easy to
	to check that
	
	(a) $\sup_{T, \theta, u} E Z_{T}^2 (\theta, u)<\infty$;
	
	(b) For all $u, v$, $\sup_{T,\theta}$ $E|Z_{T}^2(\theta,
	u)-Z_{T}^2(\theta, v)|\leq H ||u-v||^{1/2}$ for some $H>0$.
	\\
	Let $U$ be any compact closure of any open subset of $\mathbb R^l$ and
	let ${\cal Z}_{T, U} (\theta)$ and ${\cal Z}_U(\theta)$ be the same as
	${\cal Z}_{T} (\theta)$ and ${\cal Z}(\theta)$ with the integrals over ${\mathbb R^l}$ replaced by $U$.
	Then, in view of Lemma 1 and Theorem 22 of \citet{ibragimov81}, page 380, and its proof, using (a) and (b), we can have that
	for each $\theta$,
	for all $a_i\in R$
	and $\theta_i \in [0, L]$, $i=1,\ldots, k$, $k\geq 1$, as $T\rightarrow\infty$,
	\begin{eqnarray*}
		\sum_{i=1}^k a_i {\cal Z}_{T, U} (\theta_i
		)\stackrel{d}{\rightarrow}\sum_{i=1}^k a_i {\cal Z}_U(\theta_i ),
	\end{eqnarray*}
	and thus,
	\begin{eqnarray*}
		({\cal Z}_{T, U} (\theta_1), \ldots, {\cal Z}_{T, U}
		(\theta_k))^{'}\stackrel{d}{\rightarrow}({\cal Z}_U (\theta_k), \ldots,
		{\cal Z}_U (\theta_k))^{'}.
	\end{eqnarray*}
	Moreover, $\{{\cal Z}_{T, U}; T\geq 1\}$ can be viewed as a stochastic
	process in $D_{\mathbb R} ([0, L+1])$ by resetting $T+[T\theta]$, $0\leq \theta \leq L$, as $[Tu]$, $0\leq u \leq L+1$, as in \citet{billingsley68} or \citet{ethier86}.
	Using the Cauchy-Schwarz and Jensen's inequalities and the maximal moment inequality in (\ref{maximal inequality}),
	we can see that for  $r>2$ and $0\leq \theta, \theta+h \leq L$,
	\[\begin{split}
	&E\left| \int_U \left\{S_T^2(\theta+h, u)-S_T^2 (\theta,u)\right\}w(u)du \right|^r\\
	&\leq E^{1/2}\left\{\int_U |S_T(\theta+h, u)-S_T (\theta,u)|^{2r} w(u)du\right\}E^{1/2} \left\{\int_U |S_T(\theta+h, u)+S_T (\theta,u)|^{2r} w(u)du\right\}\\
	&\leq K|h|^{r/2}
	\end{split}\]
	for some $K>0$, which can yield that
	%
	$\sup_{\theta}{\cal Z}_{T, U} (\theta) \stackrel{d}{\rightarrow}
	\sup_{\theta}{\cal Z}_U(\theta)$ by virtue of Theorem 15.6 of \citet{billingsley68} and the continuous mapping
	theorem.
	Letting $U_n =\{x\in\mathbb R^l : ||x||\leq n\}$ and using the maximal moment inequality in (\ref{maximal inequality}) again, one can easily check that $\lim_{n\rightarrow\infty}\limsup_{T\rightarrow\infty}E\sup_{0\leq\theta\leq L}|{\cal Z}_T (\theta,u)-{\cal Z}_{T, U_n}(\theta, u)|=0$. Since $\sup_{\theta}{\cal Z}_{U_n}(\theta)\stackrel{d}{\rightarrow} \sup_{\theta}{\cal Z}(\theta)$ as $n\rightarrow\infty$, the theorem is validated by Proposition 6.3.9 of \citet{brockwell91}.\hfill$\Box$
	\medskip\\ \indent
	{\bf Proof of Theorem 2}. Owing to Theorem 4 of \citet{kuelbs80}, we have that for some $0<\lambda<1/2$,
	\begin{eqnarray}\label{strong conv2}\sum_{j=1}^k f(Y_j , u) -\Sigma^{1/2}(u) {\cal B}_d(k) <<
	k^{1/2-\lambda}\ \ {\rm a.s.}. \end{eqnarray} Using ($\ref{strong
		conv2}$), one can show that for each $\theta$ and $u$, as
	$T\rightarrow\infty$,
	\begin{eqnarray*}
		Z_{T}(\theta, u):=S_{T}(\theta, u)/
		q_\gamma\Big(\frac{[T\theta]}{T}\Big)\stackrel{d}{\rightarrow}Z(\theta,
		u) := \Sigma^{1/2}(u) {\cal B}_d \Big(\frac{\theta}{1+\theta}\Big)/
		\Big(\frac{\theta}{1+\theta}\Big)^\gamma.
	\end{eqnarray*}
	The rest of the proof essentially follows the same lines in that of Theorem 1 handling the univariate case
	and is omitted for brevity.  \hfill$\Box$
	

\end{document}